\documentclass[smallextended,glov3]{svjour3}
  \def\ZZ{{\mathbb Z}}
  
  \def\ZR{{\mathbb R}}

\def\beq{\begin{equation}}
\def\eeq{\end{equation}}
\def\be{\begin{equation}}
\def\ee{\end{equation}}
\def\beqar{\begin{eqnarray}}
\def\eeqar{\end{eqnarray}}
\def\ber{\begin{eqnarray}}
\def\eer{\end{eqnarray}}
\def\berb{\begin{eqnarray*}}
\def\eerb{\end{eqnarray*}}

\def\SO{\mathop{\rm SO}\nolimits}
\def\so{\mathop{\rm so}\nolimits}

\def\det{\mathop{\rm det}\nolimits}
\def\div{\mathop{\rm div}\nolimits}

\def\tr{\mathop{\rm tr}\nolimits}

\def\norm#1.#2.{\|#1\|_{#2}}
\def\Norm#1.#2.{\big\|#1\big\|_{#2}}
\def\NOrm#1.#2.{\bigg\|#1\bigg\|_{#2}}
\def\NORm#1.#2.{\Big\|#1\Big\|_{#2}}
\def\NORM#1.#2.{\Bigg\|#1\Bigg\|_{#2}}

\newcommand{\eps}{\varepsilon}

\def\vec#1{{\mathchoice{\mbox{\boldmath$\displaystyle#1$}}
{\mbox{\boldmath$\textstyle#1$}}
{\mbox{\boldmath$\scriptstyle#1$}}
{\mbox{\boldmath$\scriptscriptstyle#1$}}}}

\newcommand{\sym}{\mathop{\rm sym}\nolimits}

\def \0b{{\hbox{\boldmath $0$}}}

 \newcommand{\bbb}{{\hbox{\bf b}}}

\newcommand{\obb}{{\hbox{\bf o}}}

\newcommand{\ubb}{{\hbox{\bf u}}}

\newcommand{\ab}{\vec{a}} 
\newcommand{\cb}{\vec{c}}

\newcommand{\db}{\vec{d}}
\newcommand{\eb}{\vec{e}} \newcommand{\fb}{\vec{f}}

 \newcommand{\nb}{\vec{n}}
 
 \newcommand{\rb}{\vec{r}}
 
\newcommand{\ub}{\vec{u}} \newcommand{\vb}{\vec{v}}
 
\newcommand{\yb}{\vec{y}}

\newcommand{\Abb}{{\bf A}} \newcommand{\Bbbb}{{\bf B}}
\newcommand{\Cbb}{{\bf C}} \newcommand{\Dbb}{{\bf D}}
\newcommand{\Ebb}{{\bf E}} \newcommand{\Fbb}{{\bf F}}
\newcommand{\Gbb}{{\bf G}} \newcommand{\Hbb}{{\bf H}}
\newcommand{\Ibb}{{\bf I}}

\newcommand{\Obb}{{\bf O}} \newcommand{\Pbb}{{\bf P}}
\newcommand{\Qbb}{{\bf Q}} \newcommand{\Rbb}{{\bf R}}

\newcommand{\Ub}{\vec{U}}

\def \vphib{\vec{\varphi}} 
\def \psib{\vec{\psi}} 
 
\def \Thetab{\vec{\Theta}}

\def \Phib{\vec{\Phi}} \def \Psib{\vec{\Psi}}

\def \Thetab{\vec{\Theta}}

\newcommand{\at}{{\tilde a}}

\newcounter{primjer}[section]
\newcounter{tvrdnja}[section]
\newcounter{zadatak}[section]
\usepackage{esint}
\usepackage{amssymb}
\usepackage{graphicx}
\usepackage{enumerate}
\usepackage{verbatim}
\newcommand{\dist}{\mathop {\mbox{\rm dist}}\nolimits}

\newcommand{\meas}{\mathop {\mbox{\rm meas}}\nolimits}

\newcommand{\homeg}{\hat{\Omega}^h}

\newcommand{\ide}{\textrm{ {\it id}}}

\newcommand{\per}{\frac{x_1}{h}, \frac{x_2}{h}}

\begin{document}

\title{Periodically wrinkled plate model of the F\"oppl-von K\'arm\'an
type}




\author{
        Igor Vel\v{c}i\'{c} 
}
\institute{ Igor Vel\v{c}i\'{c} \at Faculty of Electrical Engineering and Computer Science, University of Zagreb, Unska 3, 10000 Zagreb, Croatia \\Tel: +385-1-6129965 \\Fax:+385-1-6170007 \\ \email{igor.velcic@fer.hr}}



\maketitle

\begin{abstract}
In this paper we derive, by means of $\Gamma$-convergence, the
periodically wrinkled plate model starting from three
dimensional nonlinear elasticity. We assume that the thickness
of the plate is $h^2$ and that the mid-surface of the plate is
given by $(x_1,x_2) \to (x_1,x_2,h^2\theta(\per))$, where
$\theta$ is $[0,1]^2$ periodic function. We also assume that
the strain energy of the plate has the order $h^8=(h^2)^4$,
which corresponds to the F\"oppl-von K\'arm\'an model in the
case of the ordinary plate. The obtained model mixes the
bending part of the energy with the stretching part.
\end{abstract}

\section{Introduction}

The study of thin structures is the subject of numerous works
in the theory of elasticity.  Many authors have proposed
two-dimensional shell and plate models and we come to the
problem of their justification. There is a vast literature on
the subject of plates and shells (see
\cite{Ciarlet0,Ciarlet1}).

The justification of the model of plates and shells, by using
$\Gamma$-convergence is well established. The first works in
that direction are \cite{Ledret1,LeDret2}. The thickness of the
plate is assumed to be $h$, a small parameter, and the external
loads are assumed to be of the order $0$. The obtained model
for plate and shells differs from the one obtained by the
formal asymptotic expansion in the sense that additional
relaxation of the energy functional is done.

From the pioneering work of Friesecke, James, M\"uller
\cite{Muller0} higher order models of plates and shells are
justified from three dimensional nonlinear elasticity (see
\cite{Muller0,Muller3,Muller4,Muller6,Lewicka1,Lewicka3} .
Here, higher order, relates that we assume that the magnitude
of the external loads (i.e. of the strain energy) behaves like
$h^{\alpha}, \alpha>0$ (i.e. $h^{\beta}, \beta>0$). Depending
on different parameter $\alpha$ different lower-dimensional
models are obtained (see \cite{Muller3}).

Different influence of the imperfections of the domain on the
model is also discussed in the literature. In \cite{Babadjian}
it is assumed that the stored energy function oscillates with
the order $h$, as the thickness of the plate, but the strain
energy (after divided by the order of volume $h$) is assumed to
be of the order $0$. This model thus corresponds to the one
given in \cite{Ledret1} for the ordinary plate. Also, the
influence of the specific type of the imperfections of the
domain on the F\"oppl-von K\'arm\'an plate model is discussed
in \cite{Lewicka4}. The special case of shallow shell and
weakly curved rod is discussed in \cite{Velcic1,Velcic2}. All
these models do not include periodic wrinkles which we discuss
here. Here we assume that the thickness of the plate is $h^2$
and that the mid-surface of the plate is given by $(x_1,x_2)
\to (x_1,x_2,h^2\theta(\per))$, where $\theta$ is $[0,1]^2$
periodic function. We also assume that the strain energy of the
plate (divided by the order of volume $h^2$) has the order
$h^8=(h^2)^4$, which corresponds to the F\"oppl-von K\'arm\'an
model in the case of the ordinary plate. The obtained model
mixes the bending part of the energy with the stretching part.

It could be interesting if we could generalize these periodic
wrinkles to the oscillations of the stored energy function
(like it is done in \cite{Babadjian}). But it is also important
to see that  our model depends on the pre-deformation
$\theta_0$ (see (\ref{defil0})), thus not only on the
derivatives of $\theta$. To deal with periodic wrinkles we use
the tool of two-scale convergence. The wrinkled plates model,
derived from two dimensional linear Koiter shell model, are
derived in \cite{Aganovic1,Aganovic2}. This model (and its
linearization) is different from those ones, which is expected,
since we derive the model from three dimensional nonlinear
theory and thickness of the plate is of the same order as the
amplitude of the mid-surface.

Throughout the paper  $\bar{A}$ or $\{A\}^-$  denotes the
closure of the set. By a domain we call a bounded open set with
Lipschitz boundary. $\Ibb$  denotes the identity matrix, by
$\SO(3)$ we denote the rotations in $\ZR^3$ and by $\so(3)$ the
set of antisymmetric matrices $3 \times 3$. By $\ZR^{n \times
n}_{\sym}$ we denote the set of symmetric matrices of the
dimension $n \times n$. $x'$ stands for $(x_1,x_2)$.
$\eb_1,\eb_2,\eb_3$ are the vectors of the canonical base in
$\ZR^3$. By $\ide$ we denote the identity mapping $\ide (x)=x$.
$\rightarrow$ denotes the strong convergence and
$\rightharpoonup$ the weak convergence. By $\Abb \cdot \Bbbb$
we denote $\tr(\Abb^T \Bbbb)$.
 We suppose that the
Greek indices $\alpha,\beta$ take the values in the set
$\{1,2\}$ while the Latin indices $i,j$ take the values in the
set $\{1,2,3\}$.

\section{Setting up the problem}
\setcounter{equation}{0} Let $\omega$ be a two-dimensional
domain with Lipschitz boundary in the plane spanned by $\eb_1,
\eb_2$; the generic point in $\omega$ we denote by
$x'=(x_1,x_2)$. The canonical cell in $\ZR^2$ we denote by
$Y=[0,1]^2$; the generic point in $Y$ is $y=(y_1,y_2)$. By a
periodically wrinkled plate we mean a shell defined in the
following way. Let $\theta: \ZR^2 \to \ZR$ be a $Y$-periodic
function of class $C^2$. We call $\theta$ the shape function.
We consider a three-dimensional elastic shell occupying in its
reference configuration the set $\{\hat{\Omega}^h\}^-$, where
$$\hat{\Omega}^h=\Thetab^h(\Omega^h), \ \Omega^h=\omega \times
(-\frac{h^2}{2},\frac{h^2}{2});$$ the mapping $\Thetab^h:
\{\hat{\Omega}^h\}^- \to \ZR^3$ is given by
$$ \Thetab^h(x^h)=(x_1,x_2, h^2\theta(\frac{x_1}{h},\frac{x_2}{h}))+x_3^h \nb^h(x_1,x_2)$$
for all $x^h=(x',x_3^h) \in \bar{\Omega}^h$, where $\nb^h$ is a
unit normal vector to the middle surface
$\Thetab^h(\overline{\omega})$ of the shell. By $\Omega$ we
denote $\Omega^1$ and by $x_3$ we denote $\frac{x_3^h}{h^2}$.
At each point of the surface $\bar{\omega}$ the vector $\nb^h$
is given by
$$ \nb^h(x_1,x_2)=(n^h(x_1,x_2))^{-1/2}(-h\partial_1 \theta(\frac{x_1}{h},\frac{x_2}{h}), -h \partial_2 \theta(\frac{x_1}{h},\frac{x_2}{h}), 1),$$
where  $$n^h(x_1,x_2)=h^2\partial_1
\theta(\frac{x_1}{h},\frac{x_2}{h})^2+h^2\partial_2
\theta(\frac{x_1}{h},\frac{x_2}{h})^2+1.$$ By inverse function
theorem it can be easily seen that for $h\leq h_0$ small enough
$\Thetab^h$ is a $C^1$ diffeomorphism (the global injectivity
can be proved by adapted compactness argument, see \cite[Thm
3.1-1]{Ciarlet1} for the ordinary shell). Let us by
$\theta_0:\ZR^2\to \ZR$ denote the function:
\begin{equation}
\theta_0=\theta-\langle \theta \rangle, \quad \langle \theta \rangle:= \int_Y \theta dy.
\end{equation}
 The following theorem
is easy to prove.
\begin{theorem} \label{izcea}
Then there exists $h_0=h_0(\theta)>0$ such that the Jacobian
matrix $\nabla \Thetab^h(x^h)$ is invertible for all $x^h \in
\bar{\Omega}^h$ and all $h \leq h_0$. Also there exists $C>0$
such that  for $h \leq h_0$ we have
\begin{equation} \label{determinanta} \det \nabla \Thetab^h=1+h^2 \delta^h (x^h), \end{equation}
and
\begin{equation} \label{nablaje}
\nabla \Thetab^h(x^h)=\Ibb-h \Cbb(x_1,x_2)-h^2 \Dbb(x_1,x_2,x_3)+h^3 \Rbb_1^h(x^h),
\end{equation}
\begin{equation} \label{inverz}
 (\nabla \Thetab^h(x^h))^{-1}=\Ibb+h \Cbb(x_1,x_2)+h^2 \Ebb(x_1,x_2,x_3) +h^3 \Rbb_2^h(x^h),
\end{equation}
\begin{equation} \label{glupaocjena}
\| (\nabla \Thetab^h) -\Ibb \|_{L^{\infty}(\Omega^h;\ZR^{3 \times 3})},\| (\nabla \Thetab^h)^{-1} -\Ibb \|_{L^{\infty}(\Omega^h;\ZR^{3 \times 3})}<Ch,
\end{equation}
where
\begin{equation} \label{defA}
\Cbb(x_1,x_2)=\left( \begin{array}{ccc} 0 & 0 & \partial_1 \theta(\frac{x_1}{h},\frac{x_2}{h}) \\ 0 & 0 & \partial_2 \theta(\frac{x_1}{h},\frac{x_2}{h})  \\ -\partial_1 \theta(\frac{x_1}{h},\frac{x_2}{h})  & -\partial_2 \theta(\frac{x_1}{h},\frac{x_2}{h})  & 0 \end{array} \right),
\end{equation}
\begin{eqnarray}
\nonumber & &\Dbb(x_1,x_2,x_3)=\\
& & \nonumber \left( \begin{array}{ccc} x_3 \partial_{11} \theta(\per) & x_3 \partial_{12} \theta(\per) & 0 \\ x_3 \partial_{12} \theta(\per) & x_3 \partial_{22} \theta(\per) & 0 \\ 0  & 0  & \frac{1}{2}\Big(\partial_1 \theta(\per)^2+\partial_2 \theta(\per)^2\Big) \end{array} \right), \\
\label{defB}
\end{eqnarray}
\begin{eqnarray}
 \nonumber \Ebb(x_1, x_2, x_3)&=&\Cbb^2(x_1,x_2)+\Dbb(x_1,x_2,x_3)\\ \label{defE} &=& (\Ebb_1(x_1,x_2,x_3), \Ebb_2(x_1,x_2,x_3), \Ebb_3(x_1,x_2,x_3)),
\end{eqnarray}
\begin{eqnarray}
\Ebb_1 (x_1,x_2,x_3)&= & \left(\begin{array}{c} -\partial_1 \theta (\per)^2+x_3
\partial_{11} \theta(\per) \\  -\partial_1 \theta(\per) \partial_2 \theta(\per)+  x_3 \partial_{12} \theta(\per) \\ 0
\end{array} \right), \\
\Ebb_2 (x_1,x_2,x_3)&= & \left(\begin{array}{c} -\partial_1 \theta(\per) \partial_2 \theta(\per)+  x_3 \partial_{12} \theta(\per)  \\ -\partial_2 \theta (\per)^2+x_3 \partial_{22} \theta(\per)   \\ 0
\end{array} \right), \\
\Ebb_3 (x_1,x_2,x_3)&= & \left(\begin{array}{c} 0  \\ 0   \\ - \frac{1}{2}\Big(\partial_1 \theta(\per)^2+\partial_2 \theta(\per)^2\Big)
\end{array} \right),
\end{eqnarray}

and $\delta^h: \bar{\Omega}^h \to \ZR, \
\Rbb_k^h:\bar{\Omega}^h \to \ZR^{3 \times 3},\ k=1,2$ are
functions which satisfy
$$ \sup_{0<h \leq h_0} \max_{x^h \in \bar{\Omega}^h} |\delta^h(x^h)| \leq C_0, \sup_{0<h \leq h_0} \max_{i,j} \max_{x^h \in \bar{\Omega}^h} |\Rbb_{k,ij}^h(x^h)| \leq C_0, \ k=1,2,$$
for some constant $C_0>0$.
\end{theorem}
\begin{prooof}
It is easy to see
\begin{eqnarray}
\nonumber \nb^h(x_1,x_2)&=&\eb_3-h\partial_1 \theta(\per)\eb_1-h\partial_2 \theta(\per)\eb_2\\ \label{normala111}& &\hspace{-4ex}-\frac{h^2}{2}\Big(\partial_1 \theta(\per)^2+\partial_2 \theta(\per)^2\Big)\eb_3+h^3\obb_1^h(x_1,x_2),
\end{eqnarray}
where $$ \sup_{0<h \leq h_0} \max_{x^h \in \bar{\Omega}^h}
|\obb_1^h(x_1,x_2)| \leq C, \sup_{0<h \leq h_0} \max_{x^h \in
\bar{\Omega}^h} |\partial_{\alpha} \obb_1^h(x_1,x_2)| \leq
\frac{C}{h},$$ for some $C>0$ and $\alpha=1,2$. From the
definition $\Thetab^h$ we conclude that
\begin{eqnarray} \nonumber
\nabla \Thetab^h(x',x_3^h)&=& \Ibb+h\partial_1 \theta(\per) \eb_3\otimes \eb_1+h\partial_2 \theta(\per) \eb_3\otimes \eb_2\\ & &\label{izrlamb}+(x_3^h \partial_1 \nb^h (\per), x_3^h \partial_2 \nb^h (\per), \nb^h).
\end{eqnarray}
The relation (\ref{nablaje}) is the direct consequence of the
relations (\ref{normala111}) and (\ref{izrlamb}). The relations
(\ref{determinanta}), (\ref{inverz}), (\ref{glupaocjena}) are
the direct consequences of the relation (\ref{nablaje}).
\end{prooof}

The starting point of our analysis is the minimization problem
for the wrinkled plate. The strain energy of the wrinkled plate
is given by
$$
K^{h}(\yb)= \int_{\hat{\Omega}^h} W(\nabla \yb(x))dx,
$$
where $W: \mathbb{M}^{3 \times 3} \to[0,+\infty]$ is the stored
energy density function.  $W$ is Borel measurable and, as in
\cite{Muller0,Muller3,Muller4}, is supposed to satisfy
\begin{enumerate}[i)]
\item $W$ is of class $C^2$ in a neighborhood of $\SO(3)$;
\item $W$ is frame-indifferent, i.e., $W(\Fbb)=W(\Rbb\Fbb)$ for every $\Fbb \in \ZR^{3 \times 3}$ and $\Rbb \in \SO(3)$;
\item $W(\Fbb) \geq C_W \dist^2(\Fbb,\SO(3))$, for some $C_{W}>0$ and all $\Fbb \in \ZR^{3 \times 3}$, $\ W(\Fbb)=0$ if $\Fbb \in \SO(3)$.
\end{enumerate}
By $Q_3:\ZR^{3 \times 3} \to \ZR$ we denote the quadratic form $Q_3(\Fbb)= D^2 W(\Ibb)(\Fbb,\Fbb)$ and by $Q_2: \ZR^{2 \times 2} \to \ZR$ the quadratic form,
\begin{equation}
Q_2 (\Gbb)= \min_{\ab \in \ZR^3} Q_3 (\Gbb+\ab \otimes \eb_3+\eb_3 \otimes \ab),
\end{equation}
obtained by minimizing over the stretches in the $x_3$ directions. Using ii) and iii) we conclude that both forms are positive semi-definite (and hence convex), equal to zero on antisymmetric matrices and depend only on the symmetric part of the variable matrix, i.e. we have
\begin{equation}
Q_3(\Gbb)=Q_3(\sym \Gbb), \quad Q_2 (\Gbb)= Q_2 (\sym \Gbb).
\end{equation}
Also, from ii) and iii), we can conclude that both forms are positive definite (and hence strictly convex) on symmetric matrices.
For the special case of isotropic elasticity we have
\begin{eqnarray}
\nonumber Q_3 (\Fbb) &=& 2 \mu | \frac{\Fbb+\Fbb^T}{2}|^2 + \lambda (\tr \Fbb)^2, \\
Q_2 (\Gbb) &=& 2 \mu | \frac{\Gbb+\Gbb^T}{2} |^2+ \frac{2 \mu \lambda}{2 \mu + \lambda} (\tr \Gbb)^2.
\end{eqnarray}
Since the strain energy is the most difficult part to deal with
(see Remark.....) we shall look for the $\Gamma$-limit of the
functional
$$ I^h(\yb)=\frac{1}{h^8} \frac{1}{h^2} K^h(\yb). $$
The reason why we divide by $h^8$ is because we are interested
in the F\"oppl-von K\'arm\'an type model of the wrinkled plate
(the thickness of the plate is $h^2$). The reason why we
additionally divide $K^h$ by $h^2$ is because the volume of
$\hat{\Omega}^h$ is decreasing with the order $h^2$. In the
third section we prove some technical results about two scale
convergence which we need later, in the forth section we prove
the $\Gamma$-convergence result. To prove it we firstly need
the compactness result which tells us how the displacements of
the energy order $h^8$ look like and secondly we have to prove
lower and upper bound which is standard in $\Gamma$-analysis.
\section{Two-scale convergence}
For the notion of two scale convergence see
\cite{Allaire,Nguetseng}. Here $\omega \subset \ZR^n$ is a
bounded Lipschitz domain and $Y=[0,1]^n$. For
$k=(k_i)_{i=1,\dots,n} \in \ZZ^n$ we denote by
$|k|=(\sum_{i=1}^n k_i^2)^{1/2}$.

 We denote by $C^k_{\#}(Y)$ the
space of $k$-differentiable functions with continuous $k$-th
derivative in $\ZR^2$ which are periodic of period $Y$. Then
$L^2_{\#}(Y)$ (respectively $H^m_{\#}(Y)$  is the completion
for the norm of $L^2(Y)$ (respectively $H^m(Y)$) of
$C^\infty_{\#}(Y)$. Remark that $L^2_{\#}(Y)$ actually
coincides with the space of functions in $L^2(Y)$ extended by
$Y$-periodicity to the whole of $\ZR^n$. $H^1_{\#}(Y)$
coincides with the space of functions in $H^1(Y)$ which are
$Y$-periodic at the boundary in the sense of traces and
$H^2_{\#}(Y)$  coincides with the space of functions in
$H^2(Y)$ which are with their first derivatives $Y$-periodic at
the boundary in the sense of traces etc. Using the Fourier
transform on torus it can be seen that (see e.g. \cite{Temam})
\begin{equation} \label{Fniz}
H^m_{\#} (Y)=\{u, \ u(y)=\sum_{k\in \ZZ^n} c_k e^{2 \pi i k \cdot y}, \overline{c}_k=c_{-k}, \sum_{k \in \ZZ^n} |k|^{2m} |c_k|^2 <\infty \},
\end{equation}
and the norm $\|u\|_{H^m_{\#} (Y)}$ is equivalent to the norm
$\{\sum_{k \in \ZZ^2}(1+|k|^{2m})|c_k|^2\}^{1/2}$. We also set
\begin{eqnarray} \nonumber
\dot{H}^m_{\#} (Y)=\{ u \ \textrm{of type} \ (\ref{Fniz}), \ c_0=0 \}=\{u, u \in H^m_{\#}(Y), \ \int_Y u=0 \}.
\end{eqnarray}
The restricted norm $\|u\|_{\dot{H}^m_{\#} (Y)}$ is equivalent
to the norm $\{\sum_{k \in \ZZ^2}|k|^{2m}|c_k|^2\}^{1/2}$. The
space $L^2(\omega; C_{\#}(Y))$ denotes the space of measurable
and square integrable in $x \in \omega$ with values in the
Banach space of continuous functions, $Y$-periodic in $y$. In
the analogous way one can define $L^2(\omega;H^m_{\#}(Y))$.
Using the representation (\ref{Fniz}) it can be seen that
\begin{eqnarray} \nonumber
L^2(\omega;H^m_{\#} (Y))&=&\{u, \ u(x,y)=\sum_{k\in \ZZ^n} c_k(x) e^{2 \pi i k \cdot y}, c_k \in L^2(\omega;\ZR^2), \overline{c}_k=c_{-k},\\ & & \hspace{10ex}\sum_{k \in \ZZ^n} |k|^{2m} \|c_k\|_{L^2(\omega)}^2 <\infty \},
\end{eqnarray}
where we have by $\overline{c}_k$ denoted the conjugate of
$c_k$. The norm $\|u\|_{L^2(\omega;H^m_{\#} (Y))}$ is
equivalent to the norm
$$\{\sum_{k \in
\ZZ^2}(1+|k|^{2m})\|c_k\|_{L^2(\omega)}^2\}^{1/2}.$$

The space $\mathcal{D}(\omega;C^{\infty}_{\#}(Y))$ denotes the
space of infinitely differentiable functions which take values
in $C^{\infty}_{\#}(Y)$ with compact support in $\omega$. It is
easily seen that this space is dense in $L^2(\omega;H^m_{\#}
(Y))$. In fact  it can be seen that the space of finite linear
combinations
\begin{eqnarray*}
FL(\omega; C^{\infty}_{\#}(Y))&=&\{u, \exists n \in \mathbf{N} \ u(x,y)=\sum_{k\in \ZZ^n, \ |k| \leq n} c_k(x) e^{2 \pi i k \cdot y}, \\ & & \hspace{10ex}c_k \in C_0^\infty(\omega;\ZR^2),
\overline{c}_k=c_{-k}
\}
\end{eqnarray*}
is dense in $L^2(\omega;H^m_{\#}(Y))$.
\begin{definition} \label{definicijats}
A sequence $(u_h)_{h>0}$ of functions in $L^2 (\omega)$
converges two-scale to a function $u_0$ belonging to
$L^2(\omega \times Y)$ if for every $\psi \in
L^2(\omega;C_{\#}(Y))$,
$$ \int_{\omega} u_h(x)\psi(x,\frac{x}{h}) \to \int_{\omega}
\int_Y u_0(x,y) \psi(x,y) \quad \textrm{as } h \to 0.$$
\end{definition}
By $\rightharpoonup \rightharpoonup$ we denote the two-scale
convergence. The following theorems are given in
\cite{Allaire}.
\begin{theorem} \label{L1est}
Let $f \in L^1(\omega ; C_{\#}(Y ))$. Then $f(x,\frac{x}{h})$
is a measurable function on $\omega$ for which it is valid:
\begin{equation}
\| f(x,\frac{x}{h}) \|_{L^1(\omega)} \leq \int_{\omega} \sup_{y \in Y} |f(x,y)| dy=:\|f \|_{L^1(\omega ; C_{\#}(Y ))},
\end{equation}
and
\begin{equation}
\lim_{h \to 0} \int_{\omega} f(x,\frac{x}{h}) dx=\int_\omega \int_Y f(x,y)dydx.
\end{equation}
\end{theorem}

\begin{theorem} \label{kompats}
From each bounded sequence $(u_h)_{h>0}$ in $L^2(\omega)$ one
can extract a subsequence, and there exists a limit $u_0(x,y)
\in L^2(\omega \times Y)$ such that this subsequence two-scale
converges to $u_0$.
\end{theorem}
The following theorem tells us about the form of oscillations
of order $h$ of weakly convergent sequences in $H^1(\omega)$.
\begin{theorem} \label{kompts}
Let $(u_h)$ be a bounded sequence in $H^1(\omega)$ which
converges weakly to a limit $u \in H^1(\omega)$. Then $u_h$
two-scale converges to $u(x)$, and there exists a unique
function $u_1(x,y) \in L^2(\omega; \dot{H}^1_{\#}(Y))$ such
that, up to a subsequence, $(\nabla u_h)_{h>0}$ two-scale
converges to $\nabla_x u+\nabla_y u_1(x,y)$.
\end{theorem}
\begin{remark} \label{glatkots}
In the definition of two-scale convergence we have taken the
test functions to be in the space $L^2(\omega;C_{\#}(Y))$. When
we are dealing with the sequence of the functions which are
bounded in $L^2(\omega)$ it is enough to take the test function
to be in the space $\mathcal{D}(\omega;C^{\infty}_{\#}(Y))$.
\end{remark}
The following lemmas will be needed later.
\begin{lemma}\label{prvaaaa}
Let $(u_h)_{h>0}$ be a bounded sequence in $W^{1,2}(\omega)$
and let there exists a constant $C>0$ such that $\|u_h\|_{L^2
(\omega)} \leq Ch^2$. Then we have that $(u_h)_{h>0}$ and
$(\nabla u_h)_{h>0}$ two-scale converge to 0.
\end{lemma}
\begin{prooof}
That $(u_h)_{h>0}$ two-scale converges to zero is the direct
consequence of the fact that strong convergence implies
two-scale convergence to the same limit (not depending on $y
\in Y$).  Let us now take $\psi \in
\mathcal{D}(\omega;C^{\infty}_{\#}(Y))$. Then we have
\begin{eqnarray}
\nonumber \int_{\omega} \partial_i u_h(x) \psi(x,\frac{x}{h})dx&=&-\int_{\omega} u_h(x)\partial_x \psi(x,\frac{x}{h})dx-\frac{1}{h}\int_{\omega}u_h(x) \partial_y \psi(x,\frac{x}{h})dx \\ & & \label{tsh2}
\end{eqnarray}
Since the both terms in (\ref{tsh2}) converge to $0$, due to
the fact that $\|u_h\|_{L^2(\omega)} \leq Ch^2$ we have the
claim.
\end{prooof}

The following characterization of the potentials is needed
\begin{lemma} \label{pomoccc}
Let $\ub \in L^2(\omega\times Y;\ZR^n)$ be such that for each
$\psib \in \mathcal{D}(\omega;C^{\infty}_{\#}(Y))^n$ which
satisfies $\div_y \psib=0, \ \forall x,y$ we have that
$$\int_{\Omega} \int_{Y} \ub(x,y) \psib(x,y)dydx=0.$$
Then there exists a unique function $v \in L^2(\omega;
\dot{H}^1_{\#}(Y))$ such that $\nabla_y v=u$. In the same way,
let $\Ub \in L^2(\omega\times Y;\ZR^{n \times n})$ be such that
for each $\Psib \in
\mathcal{D}(\omega;C^{\infty}_{\#}(Y))^{n\times n}$ which
satisfies $\sum_{i,j=1}^n \partial_{y_iy_j} \Psib^{ij}=0, \
\forall x,y$ we have that \begin{equation}
\label{uvjetgus}\int_{\Omega} \int_{Y} \Ub(x,y) \cdot
\Psib(x,y)dydx=0. \end{equation} Then there exists a unique
function $v \in L^2(\omega; \dot{H}^1_{\#}(Y))$ such that
$\nabla^2_y v=\Ub$.
\end{lemma}
\begin{prooof}
We shall only prove the second claim since the first goes in
the analogous way. Let us define the operator
$\nabla_y^2:L^2(\omega;\dot{H}^2_{\#}(Y)) \to L^2(\omega\times
Y;\ZR^{n \times n})$ by $v \to \nabla^2_y v$. Let us identify
the space $L^2(\omega;\dot{H}^2_{\#}(Y))$ with the sequences of
functions
\begin{eqnarray*}L^2(\omega;\dot{H}^2_{\#}(Y))&=&\{(c_k)_{k \in
\ZZ^n}, c_k \in L^2(\omega;\ZR^2), \overline{c}_k=c_{-k},
\\ & &\hspace{10ex} \int_{\omega}\sum_{k\in \ZZ^n} |k|^4
|c_k(x)|^2 dx <\infty\},
\end{eqnarray*} with the norm
$$\|(c_k)_k\|^2= \int_{\omega} \sum_{k\in
\ZZ^n} |k|^4 |c_k(x)|^2 dx $$ and the space $L^2(\omega\times
Y;\ZR^{n \times n})$ with \begin{eqnarray*}L^2(\omega \times
Y;\ZR^{n \times n}) &=& \{(c_k^{ij})_{k \in \ZZ^n, i,j=1,
\dots, n}, c_k \in L^2(\omega;\ZR^2),
\overline{c}_k^{ij}=c_{-k}^{ij},
\\ & & \hspace{10ex}\int_{\omega}\sum_{k\in \ZZ^n}
\sum_{i,j=1,\dots n}|c_k^{ij}(x)|^2 dx <\infty\},
\end{eqnarray*}
with the norm
$$\|(c_k^{ij})_k\|^2= \int_{\omega} \sum_{k\in
\ZZ^n} \sum_{i,j=1,\dots, n} |c_k^{ij} (x)|^2 dx $$
 The operator $\nabla^2_y$ operates in the
following way $$\nabla^2_y(c_k)_k= ((k_i
k_jc_k)_{i,j=1,\dots,n})_{k \in \ZZ^n}.$$ It is easily seen
that $\nabla_y^2$ is continuous and one to one. We shall prove
that it is enough to demand the condition (\ref{uvjetgus}) for
$\Psib \in FL(\omega;C^{\infty}_{\#}(Y))^{n\times n}$. Using
the properties of the Fourier transform the condition
(\ref{uvjetgus}) can be interpreted in the following way: For
given $((\Ub_k^{ij})_{i,j=1,\dots,n})_{k \in \ZZ^n} \in
L^2(\omega \times Y; (\ZR^2)^{n \times n})$ and every
$((d_k^{ij})_{i,j=1,\dots,n})_{k \in \ZZ^2} \in
FL(\omega;C^{\infty}_{\#}(Y))^{n\times n}$ which satisfies the
property \begin{equation} \sum_{i,j=1, \dots, n} k_i k_j
d_k^{ij}=0, \ \forall k \in \ZZ^n,
\end{equation}
we have that $$\int_{\omega} \sum_{k \in \ZZ^2}
\sum_{i,j=1,\dots n} \Ub_k^{ij}(x) \overline{d}_k^{ij}(x)
dx=0.$$ By fixing $k^0 \in \ZZ^n$ and taking $d_{k^0}^{ij} \in
C_0^\infty(\omega;\ZR^2)$, which satisfies $$\sum_{i,j=1,
\dots, n} k^0_i k^0_j d_{k^0}^{ij}=0$$ and defining
$d_{-k^0}^{ij}=\overline{d}_{k^0}^{ij}$, $d_k^{ij} =0,\ \forall
k \neq k^0,-k^0$, $i,j=1,\dots, n$ we conclude that
$$ \int_{\omega}\mathop{\mbox{Re}} (\Ub_{k^0}^{ij}(x)
\overline{d}_{k^0}^{ij}(x)) dx=0.$$ From this it can be easily
seen that there exists $v_{k^0} \in L^2(\omega;\ZR^2)$ such
that $\Ub_{k^0}^{ij}=k^0_i k^0_j v_{k^0}$, for all
$i,j=1,\dots, n$.   This is valid for an arbitrary $k_0 \in
\ZZ^n$ and we can easily conclude from the fact
$((\Ub_k^{ij})_{i,j=1,\dots,n})_{k \in \ZZ^2} \in L^2(\omega
\times Y; (\ZR^2)^{n \times n})$  that $\int_{\omega}\sum_{k
\in \ZZ^n} |k|^4 |v_k(x)|^2dx< \infty$. Now we have the claim
by taking $v(x,y)=(v_k(x))_k\equiv \sum_{k \in \ZZ^n}
v_k(x)e^{2 \pi i k \cdot y}$.
\end{prooof}
\begin{lemma} \label{pomoooc}
Let $(u_h)_{h>0}$ be a sequence which converges strongly
 to $u$ in $W^{1,2}(\omega)$. Let $(\vb_h)_{h>0}$ be a sequence which
 is bounded in $W^{1,2}(\omega;\ZR^n)$ and for which is valid
 \begin{equation} \label{oggrada}
 \| \nabla u_h-\vb_h\|_{L^2(\omega;\ZR^n)} \leq Ch^2,
\end{equation}
for some $C>0$. Then there exists a unique $v \in
L^2(\omega;\dot{H}^2_{\#}(Y))$ such that $\nabla \vb_h
\rightharpoonup \rightharpoonup \nabla^2 u(x)+\nabla^2_y
v(x,y)$.
\end{lemma}
\begin{prooof}
It is easily seen that $\vb_h \to \nabla u$ weakly in
$W^{1,2}(\omega; \ZR^n)$ and thus $u \in W^{2,2}(\omega)$. By
using Theorem \ref{kompats} we conclude that there exists
$\Phib^{ij} \in L^2(\omega \times Y)$ such that
$$ \nabla \vb_h^i \rightharpoonup \rightharpoonup
(\partial_{i1}u(x)+\Phib^{i1}(x,y),\dots, \partial_{in} u(x) +\Phib^{in}(x,y)).
$$ To show the existence of $v$ we shall use lemma \ref{pomoccc}.
Let us take $\Psib \in$ $\mathcal{D}(\omega;$
$C^{\infty}_{\#}(Y))^{n\times n}$ which satisfies
\begin{equation} \label{uvjetpsiii}
\sum_{i,j=1}^n
\partial_{y_iy_j} \Psib^{ij}=0, \ \forall x,y
\end{equation}
and let us calculate
\begin{eqnarray*}
\nonumber \int_{\omega} \int_{Y} \Phib \cdot  \Psib dy dx &=& \lim_{h \to 0} \int_{\omega}(\sum_{i,j=1,\dots,n}
(\partial_j \vb_h^i (x)-\partial_{ij} u(x))\cdot \Psib^{ij}(x, \frac{x}{h})) dx \\
\nonumber &=&-\lim_{h \to 0} \int_{\omega}\sum_{i,j=1,\dots,n}(\vb_h^i(x)-\partial_i u(x))\partial_{x_j} \Psib^{ij}(x,\frac{x}{h})dx \\
\nonumber& &-\lim_{h \to 0}\frac{1}{h} \int_{\omega}\sum_{i,j=1,\dots,n}(\vb_h^i(x)-\partial_i u(x))\partial_{y_j} \Psib^{ij}(x,\frac{x}{h})dx \\
\textrm{using (\ref{oggrada})} &=&
-\lim_{h \to 0}\frac{1}{h} \int_{\omega}\sum_{i,j=1,\dots,n}(\partial_i u^h (x)-\partial_i u(x))\partial_{y_j} \Psib^{ij}(x,\frac{x}{h})dx \\
&=&-\lim_{h \to 0}\frac{1}{h} \int_{\omega}\sum_{i,j=1,\dots,n}(u^h (x)-u(x))\partial_{y_jx_i} \Psib^{ij}(x,\frac{x}{h})dx \\
& &-\lim_{h \to 0}\frac{1}{h^2} \int_{\omega}\sum_{i,j=1,\dots,n}( u^h (x)- u(x))\partial_{y_i y_j} \Psib^{ij}(x,\frac{x}{h})dx\\
\textrm{using (\ref{uvjetpsiii})} &=& \lim_{h \to 0} \int_{\omega}\sum_{i,j=1,\dots,n}(\partial_j u^h (x)-\partial_j u(x))\partial_{x_i} \Psib^{ij}(x,\frac{x}{h})dx\\
& &+\lim_{h \to 0} \int_{\omega}\sum_{i,j=1,\dots,n}(u^h (x)-u(x))\partial_{x_i x_j} \Psib^{ij}(x,\frac{x}{h})dx=0.
\end{eqnarray*}
\end{prooof}
\begin{remark}
In the special case $C=0$ lemma \ref{pomoooc} is just the
generalization of  Theorem \ref{kompts}. In fact what lemma
tells us is that the claim is also valid  if we are closer to
the gradient than the order of the oscillations.
\end{remark}
\begin{lemma} \label{konvts}
Let $Q:\ZR^n \to \ZR$ be a convex function which satisfies
\begin{equation} \label{kvogrrr}
|Q(x)| \leq C(1+|x|^2), \quad \forall x\in \ZR^n,
\end{equation}
 for some $C>0$. Let $(\ub^h)_{h>0} \subset
L^2(\omega;\ZR^n)$ be a sequence which two-scale converges to
$\ub_0 \in L^2(\omega \times Y;\ZR^n)$. Then we have that
\begin{equation} \label{ocjdvk}
\int_{\omega} \int_Y Q(\ub_0(x,y))dydx \leq \liminf_{h \to 0} \int_{\omega} Q(\ub^h(x)) dx
\end{equation}
\end{lemma}
\begin{prooof}
Let us take an arbitrary $\psib \in (L^2(\omega;C_{\#}(Y)))^n$.
It is well known that if a convex function is finite on an open
set than it is continuous. Thus $Q$ is continuous. Also an
arbitrary convex function is a pointwise limit of an increasing
family of smooth convex Lipschitz functions $Q_n$. To see this
first we use the fact that there exists an increasing family
$\widetilde{Q}_n$ of piecewise affine functions (with finitely
many cuts) which pointwise converge to $Q$. Then we define
$\widetilde{\widetilde{Q_n}}=\widetilde{Q}_n-\frac{1}{n}$.
Finally we smooth every $\widetilde{\widetilde{Q}}_n$ by an
appropriate mollifier to preserve the fact that the sequence
should be increasing. We obtain for each $n \in \mathbf{N}$,
$h>0$
\begin{equation}
\int_{\omega} Q_n(\ub^h(x))dx \geq \int_{\omega}Q_n(\psib(x,\frac{x}{h}))dx+\int_{\omega} DQ_n(\psib(x,\frac{x}{h}))(\ub^h(x)-\psib(x,\frac{x}{h}))dx.
\end{equation}
By letting $h \to 0$ and using the definition of two-scale
convergence and Theorem \ref{L1est} and the fact that the
convexity of $Q_n$  and $Q_n(x) \leq C(1+|x|^2)$ implies
$|DQ_n(x)| \leq C(1+|x|)$ we obtain for each $n\in \mathbf{N}$
\begin{eqnarray} \nonumber
\liminf_{h \to 0} \int_{\omega} Q_n(\ub^h(x))dx &\geq& \int_{\omega}\int_{Y} Q_n(\psib(x,y))dy dx\\ \nonumber & &+\int_{\omega} \int_{Y} DQ_n(\psib(x,y))(\ub_0(x,y)-\psib(x,y))dx. \\ & &
\end{eqnarray}
By using an arbitrariness of $\psib$ and the density of
$L^2(\omega;C_{\#}(Y))$ in $L^2(\omega \times Y)$ we conclude
that for each  $n\in \mathbf{N}$
\begin{equation}
\liminf_{h \to 0} \int_{\omega} Q_n(\ub^h(x))dx \geq \int_{\omega} \int_Y Q_n(\ub_0(x,y))dydx.
\end{equation}
Since $Q_n<Q$ we conclude
\begin{equation}
\liminf_{h \to 0} \int_{\omega} Q(\ub^h(x))dx \geq \int_{\omega} \int_Y Q_n(\ub_0(x,y))dydx.
\end{equation}
Letting $n \to \infty$ and using (\ref{kvogrrr}) we conclude
(\ref{ocjdvk}).
\end{prooof}
\begin{lemma} \label{mnozenjets}
Let $(u_h)_{h>0}$ be a bounded sequence in $L^2(\omega)$ which
two-scale converges to $u_0(x,y) \in L^2(\omega \times Y)$. Let
$(v_h)_{h>0}$ be a sequence bounded in $L^\infty(\omega)$ which
converges in measure to $v_0 \in L^{\infty}(\omega)$. Then $v_h
u_h \rightharpoonup \rightharpoonup v_0(x) u_0(x,y)$.
\begin{prooof}
We know that $(v_h u_h)_{h>0}$ is bounded in $L^2(\omega)$ and
that there exists a subsequence of $(v_h)_{h>0}$ such that $v_h
\to v$ a.e. in $\omega$. Let us take $\psi \in \mathcal{D}
(\omega;C_{\#}(Y))$ and write,
\begin{eqnarray} \nonumber
 \int_{\omega} v_h(x)u_h
(x)\psi(x,\frac{x}{h})dx &=& \int_{\omega} (v_h(x)-v(x))u_h
(x)\psi(x,\frac{x}{h}) dx \\ & &\label{rastavvv}+\int_{\omega}
v(x)u_h (x)\psi(x,\frac{x}{h})dx
\end{eqnarray}
The second converges to
$$\int_{\omega} \int_Y v(x)u_0(x,y)\psi(x,y)dydx,$$ by the
definition of two-scale convergence. We have to prove that the
first term in (\ref{rastavvv}) converges to 0. By Egoroff's
theorem for an arbitrary $\eps>0$ there exists $E \subset
\omega$ such that $\meas(E)< \eps$ and $v_h \to v$ uniformly on
$E^c$. We write \begin{eqnarray} \nonumber \int_{\omega}
(v_h(x)-v(x))u_h (x)\psi(x,\frac{x}{h}) dx&=&\int_{E}
(v_h(x)-v(x))u_h (x)\psi(x,\frac{x}{h}) dx\\ \label{pogggg} &
&\hspace{-5ex}+\int_{E^c} (v_h(x)-v(x))u_h
(x)\psi(x,\frac{x}{h}) dx.\end{eqnarray} The second term in
(\ref{pogggg}) converges to $0$ and can be made arbitrary
small. For the first term by the Cauchy inequality we have that
there exists $C>0$ such that
\begin{equation}
\int_{E}
(v_h(x)-v(x))u_h (x)\psi(x,\frac{x}{h}) dx \leq C \sqrt{\eps} \sup_{h>0} \| u_h \|_{L^2(\omega)}.
\end{equation}
By the arbitrariness of $\eps$ we have the claim.
\end{prooof}
\end{lemma}
\section{$\Gamma$-convergence}
\setcounter{equation}{0}
We shall need the following theorem
which can be found in \cite{Muller0}.
\begin{theorem}[on geometric rigidity]\label{tgr}
Let $U \subset \ZR^m$ be a bounded Lipschitz domain, $m \geq 2$. Then there exists a constant $C(U)$ with the following property: for every $\vb\in W^{1,2}(U;\ZR^m)$ there is associated rotation $\Rbb\in \SO(m)$ such that
\begin{equation}\label{gr}
\|\nabla \vb-\Rbb\|_{L^2(U)} \leq C(U)\|\dist(\nabla \vb,\SO(m)\|_{L^2(U)}.
\end{equation}
The constant $C(U)$ can be chosen uniformly for a family of domains which are Bilipschitz equivalent with controlled Lipschitz constants.
The constant $C(U)$ is invariant under dilatations.
\end{theorem}
In the sequel we suppose $h_0 \geq \frac{1}{2}$ (see Theorem
\ref{izcea}). If this was not the case, what follows could be
easily adapted. Let us by $P^h: \Omega \to \Omega^{h} $ denote
the map $P^h(x',x_3)=(x', h^2 x_3)$. By $\nabla_h$ we denote
$$\nabla_h=\nabla_{\eb_1,\eb_2}+\frac{1}{h^2} \nabla_{\eb_3}.$$
By $\rb^h: \ZR^2 \to \ZR^2$ we denote the mapping
$\rb^h(x_1,x_2)=(\per)$. In the same way as in \cite[Theorem
10, Remark 11]{Muller3} (see also \cite[Lemma 8.1]{Lewicka1})
we can prove the following theorem. For the adaption we only
need Theorem \ref{gr} and the facts that $C(U)$ can be chosen
uniformly for Bilipschitz equivalent domains and  that the
norms $\|\nabla \Thetab^h\|$, $\|(\nabla \Thetab^h)^{-1} \|$
are uniformly bounded on $\Omega^h$ for $h \leq \frac{1}{2}$.
\begin{theorem} \label{prepisano}
Let $\omega \subset \ZR^2$ be a domain. Let $\Thetab^h$ be as
above and let  $h \leq \frac{1}{2}$. Let $\yb^h \in
W^{1,2}(\hat{\Omega}^h;\ZR^3)$ be such that
$$ \frac{1}{h^2} \int_{\homeg}\dist^2 (\nabla \yb^h,\SO(3)) dx \leq Ch^8,$$
for some $C>0$.
 Then there exists map $\Rbb^h \in
W^{1,2}(\omega,\SO(3))$ such that
\begin{eqnarray}
\| (\nabla \yb^h) \circ \Thetab^h\circ P^h -\Rbb^h \|_{L^2(\Omega)} &\leq& Ch^4, \\
 \label{ocjena1}
\| \nabla \Rbb^h \|_{L^2(\omega)} &\leq& Ch^{2}.
\end{eqnarray}
Moreover there exist a constant rotation $\bar{\Qbb}^h \in \SO(3)$ such that
\begin{eqnarray} \label{gornje}
\| (\nabla \yb^h) \circ \Thetab^h \circ P^h- \bar{\Qbb}^h\|_{L^2(\Omega)} &\leq& Ch^2,\\
 \label{pocjena}
\| \Rbb^h-\bar{\Qbb}^h \|_{L^p(\omega)} \leq C_p h^2,& & \quad \forall p<\infty.
\end{eqnarray}
Here all constants depend only on $\omega$ (and on $p$ where
indicated).
\end{theorem}
\begin{remark}
Since $\Thetab^h$ is Bilipschitz map, it can easily be seen
that the map $\yb \to \yb \circ \Thetab^h$ is an isomorphism
between the spaces $W^{1,2}(\Omega^h; \ZR^m)$ and
$W^{1,2}(\hat{\Omega}^h; \ZR^m)$ (see e.g. \cite{Adams}).
\end{remark}
To prove $\Gamma$-convergence result we need to prove the compactness result, the lower and the upper bound.
\subsection{Compactness result}
We need the following version of Korn's inequality which is proved in a standard way by contradiction.
\begin{lemma} \label{Kornova}
Let $\omega \subset \ZR^2$ be a Lipschitz domain. Then there
exists $C(\omega)>0$ such that for an arbitrary $\ub \in
W^{1,2}(\Omega; \ZR^2)$ we have
\begin{equation} \label{Korn}
\| \ub \|_{W^{1,2}(\omega;\ZR^2)} \leq C(\omega) (\|\sym \nabla \ub \|_{L^2(\omega;\ZR^2)}+|\int_{\omega} \ub dx|+|\int_{\omega} (\partial_2 \ub_{1}-\partial_1 \ub_{2})dx|).
\end{equation}
\end{lemma}
\begin{lemma} \label{osnovna}
Let $\yb^h \in W^{1,2}(\hat{\Omega}^h;\ZR^3)$ be such that
\begin{eqnarray}
\frac{1}{h^2} \int_{\homeg}\dist^2 (\nabla \yb^h,\SO(3)) dx & \leq & C h^8,
\end{eqnarray}
Then there exists maps $\Rbb^h$  $\in$  $W^{1,2}(\omega,
\SO(3))$ and constants $\bar{\Rbb}^h \in \SO(3)$, $\cb^h \in
\ZR^3$  such that
$$ \widetilde{\yb}^h:=(\bar{\Rbb}^h)^T \yb^h- \cb^h$$
and the corrected in-plane and the out-of-plane displacements
\begin{eqnarray} \nonumber \ub^h (x_1,x_2)&:=&
\frac{1}{h^4} \Bigg(\int_{-1/2}^{1/2} \left( \begin{array}{c} \widetilde{\yb}^h_1 \circ \Thetab^h \circ P^h \\  \widetilde{\yb}^h_2 \circ \Thetab^h \circ P^h \end{array} \right)(x_1,x_2,x_3)dx_3-\left(\begin{array}{c} x_1 \\x_2\end{array}\right)\\ \nonumber& &-h^2\left(\begin{array}{c}\Rbb^h_{13}(x_1,x_2) \theta_0(\per) \\ \Rbb^h_{23}(x_1,x_2)\theta_0(\per)\end{array}\right) \Bigg), \end{eqnarray} \begin{eqnarray}  \label{definicijauv} v^h(x_1,x_2)&:=&
 \frac{1}{h^2}\int_{-1/2}^{1/2} (\widetilde{\yb}_3^h \circ \Thetab^h\circ P^h)(x_1,x_2,x_3)dx_3- \theta(\per) \end{eqnarray}
satisfy
\begin{equation} \label{sranje1}
\| (\nabla \widetilde{\yb}^h) \circ \Thetab^h \circ P^h-\Rbb^h \|_{L^2(\Omega)} \leq C h^4,
\end{equation}
\begin{equation}\label{sranje2}
\| \Rbb^h-\Ibb\|_{L^p(\omega)} \leq C_p h^2 \quad \forall p<\infty, \quad \| \nabla \Rbb^h \|_{L^2(\omega)} \leq Ch^2.
\end{equation}
Moreover every subsequence (not relabeled)  has its subsequence (also not relabeled) such that
\begin{eqnarray}\label{konvergencija}
v^h &\to& v  \quad \textrm{in } W^{1,2}(\omega), \quad v \in W^{2,2}(\omega) \\ \label{konvergencija2}
\ub^h &\rightharpoonup& \ub \quad \textrm{in } W^{1,2}(\omega;\ZR^2), \\ \label{acc}
\frac{\Rbb^h-\Ibb}{h^2} &\rightharpoonup& \Abb \quad \textrm{in } W^{1,2}(\omega;\ZR^{3 \times 3}),
\\ \label{graa}
\frac{(\nabla \widetilde{\yb}^h)\circ \Thetab^h \circ P^h-\Ibb}{h^2} &\to& \Abb \quad \textrm{in } L^2(\Omega;\ZR^{3 \times 3})
\\
\partial_3 \Abb=0, &\quad& \Abb \in W^{1,2} (\omega; \ZR^{3 \times 3}),
\\ \label{defabb}
\Abb&=&\eb_3 \otimes \nabla v -\nabla v \otimes \eb_3, \\
 \label{simetricno} \frac{\sym (\Rbb^h-\Ibb)}{h^4} &\to&
\frac{\Abb^2}{2} \quad \textrm{in } L^2(\Omega;\ZR^{3 \times
3}).
\end{eqnarray}
\end{lemma}
\begin{prooof}
We shall follow the proof of Lemma 13 in \cite{Muller3} (see
also Lemma 2 in \cite{Velcic1}). Estimates (\ref{sranje1}) and
(\ref{sranje2}) follow immediately from Theorem \ref{prepisano}
since one can choose $\bar{\Rbb}^h$ so that (\ref{gornje})
holds with $\bar{\Qbb}^h=\Ibb$. By applying additional constant
in-plane rotation of order $h^2$ to $\widetilde{\yb}^h$ and
$\Rbb^h$  we may assume in addition to (\ref{sranje1}) and
(\ref{sranje2}) that
\begin{equation} \label{josmalo}
\int_{\Omega} ((\partial_2 \widetilde{\yb}^h_{1})\circ \Thetab^h \circ P^h-(\partial_1 \widetilde{\yb}^h_{2})\circ \Thetab^h \circ P^h) dx=0.
\end{equation}
By choosing $\cb^h$ suitably we may also assume that
\begin{equation} \label{nacinzac}
\int_{\Omega} (\widetilde{\yb}^h\circ \Thetab^h \circ P^h-\Thetab^h \circ P^h)dx=0.
\end{equation}
Let us define $\Abb^h=\frac{\Rbb^h-\Ibb}{h^2}$. From
(\ref{sranje2}) we get for a subsequence
$$ \Abb^h \rightharpoonup \Abb \quad \textrm{in } W^{1,2}(\omega;\ZR^{3 \times 3}) $$
Thus we deduce (\ref{acc}).  Using (\ref{glupaocjena}),
(\ref{sranje1}) we deduce (\ref{graa}). Since $(\Rbb^h)^T
\Rbb^h=\Ibb$ we have $\Abb^h+(\Abb^h)^T=-h^2(\Abb^h)^T \Abb^h$.
Hence $\Abb+\Abb^T=0$ and after multiplication by $1/h^2$ we
obtain (\ref{simetricno}) from the strong convergence of
$\Abb^h$. Using (\ref{sranje1}) we conclude that
\begin{equation} \label{sranje10}
\| \sym((\nabla \widetilde{\yb}^h) \circ \Thetab^h \circ P^h-\Ibb) \|_{L^2(\Omega)} \leq C h^4.
\end{equation}
The following is useful
\begin{eqnarray} \nonumber
(\nabla \widetilde{\yb}^h) \circ \Thetab^h \circ P^h&=&(\nabla(\widetilde{\yb}^h \circ \Thetab^h)\circ P^h) ((\nabla \Thetab^h)^{-1} \circ P^h)\\ \label{kljucnozadokaz1}&=&\nabla_h (\widetilde{\yb}^h \circ \Thetab^h \circ P^h) ((\nabla \Thetab^h)^{-1} \circ P^h).
\end{eqnarray}
From (\ref{kljucnozadokaz1}) it follows
\begin{equation} \label{kljucnozadokaz2}
((\nabla \widetilde{\yb}^h) \circ \Thetab^h \circ P^h)((\nabla \Thetab^h) \circ P^h)=\nabla_h (\widetilde{\yb}^h \circ \Thetab^h \circ P^h).
\end{equation}
Using (\ref{glupaocjena}), (\ref{graa}),
(\ref{kljucnozadokaz2}) we conclude that
\begin{equation} \label{zazazaKorna}
\frac{1}{h^2} \nabla_h (\widetilde{\yb}^h \circ \Thetab^h \circ P^h-\Thetab^h \circ P^h) \to \Abb \ \textrm{in } L^2(\Omega;\ZR^{3 \times 3}).
\end{equation}
From (\ref{nacinzac}), (\ref{zazazaKorna}) and the Poincare
inequality we conclude the convergence in
(\ref{konvergencija}).  Moreover we have $\partial_i
v=\Abb_{3i}$ for $i=1,2$. Hence $v \in W^{2,2}$ since $\Abb \in
W^{1,2}$. Since $\Abb$ is skew-symmetric we immediately have
$\Abb_{13}=-\partial_1 v, \Abb_{23}=-\partial_2 v$.

By the chain rule the following identities are valid for
$i=1,2,3$, $\alpha=1,2$,
\begin{eqnarray}
\nonumber\partial_\alpha (\widetilde{\yb}^h_i \circ \Thetab^h \circ P^h)&=& (\partial_1 \widetilde{\yb}^h_i) \circ \Thetab^h \circ P^h \cdot (\partial_\alpha \Thetab^h_1)\circ P^h\\ \nonumber & & \hspace{-23ex}+(\partial_2 \widetilde{\yb}^h_i) \circ \Thetab^h \circ P^h \cdot  (\partial_\alpha \Thetab^h_2) \circ P^h +
(\partial_3 \widetilde{\yb}_i^h) \circ \Thetab^h \circ P^h \cdot (\partial_\alpha \Thetab^h_3) \circ P^h \\ \label{parcal} & & \\
\nonumber \frac{1}{h^2} \partial_3 (\widetilde{\yb}^h_i \circ \Thetab^h \circ P^h)&=&(\partial_1 \widetilde{\yb}^h_i) \circ \Thetab^h \circ P^h \cdot (\partial_3 \Thetab^h_1)\circ P^h\\ \nonumber & & \hspace{-23ex}+(\partial_2 \widetilde{\yb}^h_i) \circ \Thetab^h \circ P^h \cdot  (\partial_3 \Thetab^h_2) \circ P^h +
(\partial_3 \widetilde{\yb}_i^h) \circ \Thetab^h \circ P^h \cdot (\partial_3 \Thetab^h_3) \circ P^h \\ \label{parc3}& &
\end{eqnarray}
From (\ref{nablaje}), (\ref{graa}), (\ref{sranje10}) we
conclude for $\alpha,\beta=1,2$
\begin{equation}
\|(\partial_\alpha \widetilde{\yb}^h_\alpha) \circ \Thetab^h \circ P^h \cdot (\partial_\beta \Thetab^h_\beta)\circ P^h-(\partial_\beta \Thetab^h_\beta)\circ P^h\|_{L^2(\Omega)} \leq Ch^4,
\end{equation}
and
\begin{equation}
\|(\partial_3 \widetilde{\yb}^h_\alpha) \circ \Thetab^h \circ P^h \cdot (\partial_\alpha \Thetab^h_3)\circ P^h-h \Rbb_{\alpha 3}(\partial_\alpha \theta) \circ \rb^h\|_{L^2(\Omega)} \leq Ch^4,
\end{equation}
 for some $C>0$.
Using (\ref{sranje1}), (\ref{sranje2}), (\ref{graa}),
(\ref{josmalo}), (\ref{parcal}) we conclude that
\begin{equation} \label{antisimmm}
|\int_{\omega} (\partial_2 \ub^h_1- \partial_1 \ub^h_2) dx | \leq Ch^4.
\end{equation}
In the same way using (\ref{sranje1}), (\ref{sranje2}),
(\ref{graa}), (\ref{josmalo}), (\ref{parcal}) we conclude that
$\| \frac{1}{h^4} \sym \nabla \ub^h \|_{L^2(\omega)}$ is
bounded. Using Lemma \ref{Kornova}, (\ref{nacinzac}) and
(\ref{antisimmm}) we have the convergence
(\ref{konvergencija2}).
 It remains to conclude $\Abb_{12}=0$. But this is easy, from  (\ref{konvergencija2}), (\ref{acc}) and (\ref{zazazaKorna}).
\end{prooof}
\begin{lemma} \label{ovavecer}
If we additionally assume the following
\begin{eqnarray} \label{uuuuu}
 \Bigg\|\int_{-1/2}^{1/2} \left( \begin{array}{c} \widetilde{\yb}^h_1 \circ \Thetab^h \circ P^h \\  \widetilde{\yb}^h_2 \circ \Thetab^h \circ P^h \end{array} \right)(x_1,x_2,x_3)dx_3-\left(\begin{array}{c} x_1 \\x_2\end{array}\right)\Bigg\|_{L^2(\omega)}&\leq& Ch^4, \\ \label{vvvvv}
 \big\|\int_{-1/2}^{1/2} (\widetilde{\yb}_3^h \circ \Thetab^h\circ P^h)(x_1,x_2,x_3)dx_3\big\|_{L^2(\omega)} &\leq& Ch^2, \end{eqnarray}
we can take $\bar{\Rbb}^h=\Ibb$, $\cb^h=0$ in Lemma
\ref{osnovna}.
\end{lemma}
\begin{prooof}
First we shall prove that we can take $\bar{\Qbb}^h=\Ibb$ in
Theorem \ref{prepisano}. By multiplying (\ref{gornje}) with
$\nabla (\Thetab^h \circ P^h)$ and using Poincare inequality on
$\Omega$ we have that there exists $ b^h \in \ZR^3$ such that

\begin{equation} \label{jjjjjj}\| \yb^h \circ \Thetab^h \circ P^h- \bar{\Qbb}^h
(\Thetab^h\circ P^h)-\bbb^h\|_{L^2(\Omega)} \leq Ch^2.
\end{equation} From (\ref{uuuuu}), (\ref{vvvvv}) and
by integrating (\ref{jjjjjj})  with respect to $x_3$ we
conclude
\begin{equation}
\left\|\bar{\Qbb}^h \left(\begin{array}{c} x_1 \\x_2 \\0 \end{array}\right) - \left(\begin{array}{c} x_1 \\x_2 \\ 0\end{array}\right)-\bbb^h \right\|_{L^2(\omega)} \leq Ch^2.
\end{equation}
From this, using the fact that $\bar{\Qbb^h} \in \SO(3)$, we
can conclude that $\| \bbb^h \| \leq Ch^2$ and
$\|\bar{\Qbb}^h-\Ibb\| \leq Ch^2$. Thus we can take
$\bar{\Qbb}^h=\Ibb$ in Theorem \ref{prepisano}. Now we conclude
that we can repeat the proof of Lemma \ref{osnovna} with
assumptions (\ref{uuuuu}) and (\ref{vvvvv}) instead of
assumptions (\ref{josmalo}) (i.e. (\ref{antisimmm})) and
(\ref{nacinzac}).

\end{prooof}

\begin{lemma} \label{kompor}
Let $\widetilde{\yb}^h, \Rbb^h, \ub^h, \ub,v^h, v$ be as in
Lemma \ref{osnovna}. Then there exist unique $\ub_1 \in
L^2(\omega ; \dot{H}^1_{\#}(Y))^2$ and $v_1 \in L^2(\omega;
\dot{H}^2_{\#}(Y))$ such that
\begin{eqnarray} \label{dvouh}
\nabla \ub^h &\rightharpoonup \rightharpoonup& \nabla_x \ub(x)+\nabla_y \ub_1 (x,y) \textrm{ in } L^2(\omega \times Y;\ZR^{2 \times 2}), \\ \label{rds}
\nabla \frac{1}{h^2} \left( \begin{array}{c} \Rbb^h_{13} \\ \Rbb^h_{23} \end{array} \right) & \rightharpoonup \rightharpoonup& -\nabla^2_x v(x)-\nabla^2_y v_1(x,y) \textrm{ in } L^2(\omega \times Y; \ZR^{2 \times 2}).
\end{eqnarray}
\end{lemma}
\begin{prooof}
The relation (\ref{dvouh}) is the direct consequence of the
relation (\ref{konvergencija2}) and Theorem \ref{kompts}. To
prove (\ref{rds}) we shall use Lemma \ref{pomoccc}. From the
relation (\ref{sranje1}) using the boundedness of
$\nabla\Thetab^h$ and (\ref{kljucnozadokaz2}) we conclude
\begin{eqnarray} \nonumber
& &\| \frac{1}{h^2}\int_{-1/2}^{1/2}(\nabla_h (\widetilde{\yb}^h \circ \Thetab^h \circ P^h-\Thetab^h \circ P^h)dx_3\\\label{zadanasss}  & &\hspace{10ex}-\frac{1}{h^2} (\Rbb^h-\Ibb)\int_{-1/2}^{1/2}\nabla_h(\Thetab^h \circ P^h)dx_3
\|_{L^2(\omega)} \leq C h^2.
\end{eqnarray}
From (\ref{zadanasss}) using (\ref{nablaje}) and
(\ref{simetricno}) we conclude
\begin{equation}
\|\partial_1 v^h-\frac{1}{h^2} \Rbb^h_{31}\|_{L^2(\omega)} \leq Ch^2, \  \|\partial_2 v^h- \frac{1}{h^2} \Rbb^h_{32}\|_{L^2(\omega)} \leq Ch^2,
\end{equation}
for some $C>0$. From that we have, by using (\ref{simetricno}),
\begin{equation}
\|\partial_1 v^h+\frac{1}{h^2} \Rbb^h_{13}\|_{L^2(\omega)} \leq
Ch^2, \ \|\partial_2 v^h+ \frac{1}{h^2}
\Rbb^h_{23}\|_{L^2(\omega)} \leq Ch^2.
\end{equation}
The claim is now the direct consequence of Lemma \ref{pomoccc}
and (\ref{acc}).

\end{prooof}
\subsection{Lower bound}
\begin{remark} \label{prosirenadk}
In the next lemma we shall characterize the limiting strain and
express it in terms of $\ub, \ub_1, v,v_1,\theta$ (see Lemma
\ref{kompor}). Since the strain is an element of $L^2(\Omega)$
and we want to obtain its two-scale limit which characterizes
only oscillations of the first two variables of the order $h$
the following modification of Definition \ref{definicijats}) is
needed: A sequence $(u_h)_{h>0}$ of functions in $L^2 (\Omega)$
converges two-scale to a function $u_0$ belonging to
$L^2(\Omega \times Y)$ if for every $\psi \in
L^2(\Omega;C_{\#}(Y))$,
$$ \int_{\Omega} u_h(x_1,x_2,x_3)\psi(x_1,x_2,x_3,\per) \to \int_{\Omega}
\int_Y u_0(x_1,x_2,x_3,y) \psi(x_1,x_2,x_3,y).$$ Here
$Y=[0,1]^2$. It can be also seen that the analogous statements
of Theorem \ref{L1est}, Theorem \ref{kompats}, Theorem
\ref{kompts} and Remark \ref{glatkots} are valid (see
\cite{Pavliotis} for time dependent problems). Also the
analogous conclusions of Lemma \ref{prvaaaa} Lemma
\ref{pomoccc}, Lemma \ref{pomoooc}, Lemma \ref{konvts} and
Lemma \ref{mnozenjets} are valid ($\nabla_x$ should be replaced
by the gradient in the first two variables).
\end{remark}
\begin{lemma} \label{identifikacija}
Consider $\yb^h:\hat{\Omega}^h \to \ZR^3$, $\Rbb^h \in
W^{1,2}(\omega; \SO(3))$ and define $\ub^h,v$ by
(\ref{definicijauv}) Suppose that we have a subsequence of
$\yb^h$ such that (\ref{sranje1})-(\ref{simetricno}) are valid.
Additionally we suppose
\begin{eqnarray} \label{cudno2}
& &\nabla \ub^h \rightharpoonup \rightharpoonup \nabla_x \ub(x)+\nabla_y \ub_1 (x,y) \textrm{ in } L^2(\omega \times Y;\ZR^{2 \times 2}), \\ \label{pprrts} & &
\nabla \frac{1}{h^2} \left( \begin{array}{c} \Rbb^h_{13} \\
\Rbb^h_{23} \end{array} \right)  \rightharpoonup
\rightharpoonup -\nabla^2_x v(x)-\nabla^2_y v_1(x,y) \textrm{
in } L^2(\omega \times Y; \ZR^{2 \times 2}),
\end{eqnarray}
for  $\ub_1 \in L^2(\omega ; \dot{H}^1_{\#}(Y))^2$ and $v_1 \in
L^2(\omega; \dot{H}^2_{\#}(Y))$. Then
\begin{equation}
\Gbb^h:= \frac{(\Rbb^h)^T ((\nabla \yb^h) \circ \Thetab^h\circ P^h)-\Ibb}{h^4} \rightharpoonup \rightharpoonup \Gbb \quad \textrm{in } L^2(\Omega \times Y;\ZR^{3 \times 3}),
\end{equation}
and the $2 \times 2$ sub-matrix $\Gbb''$ given by $\Gbb''_{\alpha \beta}=\Gbb_{\alpha \beta}$ for $1 \leq \alpha,\beta \leq 2$ satisfies
\begin{equation} \label{zakor}
\Gbb''(x',x_3)=\Gbb_0(x_1,x_2,y)+x_3 \Gbb_1(x_1,x_2,y),
\end{equation}
where
\begin{eqnarray} \nonumber
\sym \Gbb_0(x_1,x_2,y)&=& \sym\nabla_x \ub(x)+\sym\nabla_y \ub_1(x,y)\\ \label{defg0} & & \hspace{-5ex}-\nabla_x^2 v(x) \theta_0(y)-\nabla^2_y v_1(x,y)\theta_0(y)+\frac{1}{2}\nabla_x v \otimes \nabla_x v \\ \label{defg1}
\Gbb_1(x_1,x_2,y)&=&-(\nabla_x)^2 v(x)-(\nabla_y)^2 v_1(x,y).
\end{eqnarray}
\end{lemma}
\begin{prooof} We follow the proof of Lemma 15 in \cite{Muller3} (see also Lemma 4 in \cite{Velcic1}).
By the assumption $\Gbb^h$ is bounded in $L^2$, thus a
subsequence converges weakly.

To show that the limit matrix $\Gbb''$ is affine in $x_3$ we
consider the difference quotients
\begin{equation}
\Hbb^h(x',x_3)=s^{-1}[\Gbb^h(x_1,x_2,x_3+s)-\Gbb^h(x_1,x_2,x_3)].
\end{equation}
By multiplying the definition of $\Gbb^h$ with $\Rbb^h$ and using (\ref{kljucnozadokaz2}) we obtain for $\alpha,\beta \in \{1,2\}$
\begin{eqnarray} \nonumber
(\Rbb^h \Hbb^h)_{\alpha \beta} &=& \frac{1}{s h^4} [\nabla_h (\yb^h \circ  \Thetab^h \circ P^h)(x',x_3+s) -\\ \nonumber & & \hspace{10ex}  -\nabla_h (\yb^h \circ  \Thetab^h \circ P^h)(x',x_3)]_{\alpha \beta} \\  \nonumber & &
+\frac{1}{s h^4}\Big[((\nabla \yb^h) \circ \Thetab^h \circ P^h)(x',x_3+s)\cdot\\ \nonumber & & \hspace{10ex} \cdot \left(\Ibb-((\nabla \Thetab^h) \circ P^h)(x',x_3+s)\right)\Big]_{\alpha \beta} \\ \nonumber
& &-\frac{1}{s h^4}\Big[((\nabla \yb^h) \circ \Thetab^h \circ P^h)(x',x_3) \cdot \\ & & \nonumber \hspace{10ex} \cdot \left(\Ibb-((\nabla \Thetab^h) \circ P^h)(x',x_3)\right)\Big]_{\alpha \beta} \\ \nonumber
&=&  \frac{1}{s h^4} \Big[(\nabla_h (\yb^h \circ  \Thetab^h \circ P^h-\Thetab^h \circ P^h)(x',x_3+s)-\\ \nonumber & & \hspace{10ex}  -\nabla_h (\yb^h \circ  \Thetab^h \circ P^h-\Thetab^h \circ P^h)(x',x_3)\Big]_{\alpha \beta} \\ \nonumber
& &+\frac{1}{s h^4}\Big[((\nabla \yb^h) \circ \Thetab^h \circ P^h)(x',x_3+s)-\\ \nonumber & & \hspace{10ex} -\Rbb^h(x'))\left(\Ibb-((\nabla \Thetab^h) \circ P^h)(x',x_3+s)\right)\Big]_{\alpha \beta} \\ \nonumber
& &-\frac{1}{s h^4}\Big[((\nabla \yb^h) \circ \Thetab^h \circ P^h)(x',x_3)\\ \nonumber& &\hspace{10ex} -\Rbb^h(x'))\left(\Ibb-((\nabla \Thetab^h) \circ P^h)(x',x_3)\right)\Big]_{\alpha \beta}
 \\ \nonumber
& &-\frac{1}{s h^4}\Big[(\Rbb^h (x')-\Ibb) ((\nabla \Thetab^h) \circ P^h)(x',x_3+s)-\\ & & \hspace{10ex} \label{zadnjarel}-((\nabla \Thetab^h) \circ P^h)(x',x_3)\Big]_{\alpha \beta}.
\end{eqnarray}
By using (\ref{nablaje}), (\ref{sranje1}) we conclude that the
second and third term converges to $0$ strongly in $L^2(\omega
\times (-\frac{1}{2},\frac{1}{2}-s)$. The forth term also
converges to $0$ by (\ref{nablaje}) and (\ref{acc}). Thus we
have that the first term in (\ref{zadnjarel}) is bounded in
$L^2(\Omega)$ and thus two scale converges to some $\Gbb_1 \in
L^2(\Omega \times Y)$. We conclude
\begin{eqnarray}\nonumber
& &\frac{1}{s h^4} [(\nabla_h (\yb^h \circ  \Thetab^h \circ P^h-\Thetab^h \circ P^h)(x',x_3+s)\\ \nonumber & & \hspace{15ex}-\nabla_h (\yb^h \circ  \Thetab^h \circ P^h-\Thetab^h \circ P^h)(x',x_3))]_{\alpha \beta}\\& & \hspace{7ex} =\frac{1}{h^2} \partial_{\beta} \Big( \frac{1}{s} \int_0^s \frac{1}{h^2} \partial_3 (\yb^h \circ  \Thetab^h \circ P^h-\Thetab^h \circ P^h)_{\alpha} \Big).
\end{eqnarray}
By using (\ref{sranje10}) (which is a consequence of
(\ref{sranje1}) and (\ref{simetricno})) and (\ref{parc3}) we
conclude that there exists $C>0$ such that for $\alpha=1,2$,
$\beta=3-\alpha$
\begin{eqnarray} \nonumber
& &\Big\| \Big( \frac{1}{s h^2} \int_0^s \frac{1}{h^2} \partial_3 (\yb^h \circ  \Thetab^h \circ P^h-\Thetab^h \circ P^h)_{\alpha} \Big)-\frac{1}{h^2} \Rbb^h_{\alpha 3}\\ & &\hspace{10ex} +\frac{1}{h}\Rbb^h_{\alpha \beta}(\partial_{\beta}\theta)\circ \rb^h \Big\|_{L^2(\omega\times (-\frac{1}{2},\frac{1}{2}-s))} <Ch^2.
\end{eqnarray}
By using Lemma \ref{prvaaaa}, the fact that for $\alpha=1,2$,
$\beta=3-\alpha$
\begin{equation}
\frac{1}{h}\Rbb^h_{\alpha \beta}(\partial_{\beta}\theta) \circ \rb^h \to 0 \textrm{ in } W^{1,2}(\omega)
\end{equation}
(this follows from (\ref{acc})),
 we conclude from
(\ref{zadnjarel}) that
\begin{equation}
\Rbb^h \Hbb^h \rightharpoonup \rightharpoonup \Gbb_1, \quad \textrm{in } L^2(\omega \times (-\frac{1}{2}, \frac{1}{2}-s);\ZR^{n \times n}),
\end{equation}
where $\Gbb_1$ is given by (\ref{defg1}). Since $\Rbb^h \to
\Ibb$ boundedly a.e. we conclude $\Hbb^h \rightharpoonup
\rightharpoonup \Gbb_1$ from Lemma \ref{mnozenjets}. From this
we have also (\ref{zakor}). In order to prove formula for
$\Gbb_0$ it suffices to study
$$ \Gbb_0^h (x')=\int_{-\frac{1}{2}}^{\frac{1}{2}} \Gbb^h(x',x_3) dx_3. $$
We have for $\alpha, \beta \in \{1,2\}$
\begin{eqnarray} \nonumber
(\Gbb^h)_{\alpha \beta}(x',x_3) &=&  \frac{((\nabla \yb^h)\circ \Thetab^h \circ P^h-\Ibb)_{\alpha \beta}}{h^4}-\frac{(\Rbb^h-\Ibb)_{\alpha \beta}}{h^4}\\  \label{idensim}& &+\Big[(\Rbb^h-\Ibb)^T \frac{(\nabla \yb^h) \circ \Thetab^h \circ P^h-\Rbb^h}{h^4} \Big]_{\alpha \beta}.
\end{eqnarray}
The third term in (\ref{idensim}) converges strongly to $0$ in
$L^2(\Omega)$. From (\ref{parcal}) we conclude, after a little
calculation by using (\ref{nablaje}), (\ref{sranje1}) and
(\ref{acc}),that
\begin{eqnarray}\nonumber
& & \nabla \ub^h-  \int_{-1/2}^{1/2} \frac{((\nabla \yb^h)\circ \Thetab^h \circ P^h-\Ibb)}{h^4}dx_3 \\ \label{dodatak1}
& &+(\theta_0 \circ \rb^h) \nabla \frac{1}{h^2} \left( \begin{array}{c} \Rbb^h_{13} \\
\Rbb^h_{23} \end{array} \right) \to 0 \textrm{ in } L^2(\Omega;\ZR^{2 \times 2}).
\end{eqnarray}
From (\ref{dodatak1}) we conclude
\begin{eqnarray} \nonumber
\int_{-1/2}^{1/2} \frac{((\nabla \yb^h)\circ \Thetab^h \circ P^h-\Ibb)}{h^4}dx_3 &\rightharpoonup \rightharpoonup&  \nabla_x \ub(x)+\nabla_y \ub_1(x,y)\\& &\label{dodaaaa} \hspace{-33ex}-(\nabla_x)^2 v(x)\theta_0(y)-(\nabla_y)^2 v_1(x,y)\theta_0(y) \ \textrm{ in } L^2(\omega \times Y;\ZR^{2 \times 2}).
\end{eqnarray}
From (\ref{dodatak1}) and (\ref{dodaaaa}) we have
(\ref{defg0}).
\end{prooof}

By $I^L:W^{1,2}(\omega;\ZR^2) \times
L^2(\omega;\dot{H}_{\#}^1(Y))^2 \times W^{2,2}(\omega) \times
L^2(\omega;\dot{H}_{\#}^2) \to \ZR_0^+ $ we denote the
functional
\begin{eqnarray} \nonumber
I^{L}(\ub,\ub_1,v, v_1)&=&\int_{\omega} \int_Y \Big( \frac{1}{2} Q_2\big(\sym\nabla_x \ub(x)+\sym\nabla_y \ub_1(x,y)-\nabla_x^2 v(x) \theta_0(y)\\ \nonumber& & \hspace{10ex}-\nabla^2_y v_1(x,y)\theta_0(y)+\frac{1}{2}\nabla_x v \otimes \nabla_x v\big)\Big)dydx\\ \label{defL}& &\hspace{-10ex} +\frac{1}{24}\int_\omega \int_Y Q_2(\nabla_y^2 v_1(x,y))dy dx
+ \frac{1}{24}\int_\omega Q_2(\nabla_x^2 v(x)) dx.
\end{eqnarray}

\begin{corollary} \label{zajakuisto}
Let $\yb^h,\Rbb^h,\ub^h,\ub,\ub_1,v^h,v,v_1, \Gbb^h, \Gbb,
\Gbb'', \Gbb_0, \Gbb_1$ be as in  Lemma \ref{identifikacija}.
Then we have the following semi-continuity results.
\begin{equation}
\liminf_{h \to 0} \frac{1}{h^8} I^h (\yb^h) \geq I^L(\ub,\ub_1,v,v_1).
\end{equation}
\end{corollary}
\begin{prooof}
 We shall use the truncation, the Taylor expansion and the weak semi-continuity argument as in the proof of Corollary $16$ in \cite{Muller3}.
Let $m:[0,\infty) \to [0,\infty)$ denote a modulus of
continuity of $D^2W$ near the identity and consider the good
set $\Omega_h :=\{x \in \Omega: |\Gbb^h (x)| < h^{-1} \}$. Its
characteristic function $\chi_h $ is bounded and satisfies
$\chi_h \to 1$ a.e. in $\Omega$. Thus we have $\chi_h \Gbb^h
\rightharpoonup \rightharpoonup \Gbb$ in $L^2 (\Omega \times
Y;\ZR^{3 \times 3})$ by Lemma \ref{mnozenjets}. By Taylor
expansion
\begin{equation} \label{uvjetnaniz}
\frac{1}{h^8} \chi_h W(\Ibb+h^4 \Gbb^h) \geq \frac{1}{2} Q_3(\chi_h \Gbb^h)-m(h^3)|\Gbb^h|^2.
\end{equation}
Using i), ii), (\ref{uvjetnaniz}), the boundedness of sequence
$\Gbb^h$ in $L^2(\Omega;\ZR^{3 \times 3})$ and Lemma
\ref{konvts} we conclude
\begin{eqnarray}
\nonumber & & \liminf_{h \to 0} \frac{1}{h^8} I^h(\yb^h)
= \liminf_{h \to 0} \frac{1}{h^8} \int_{\Omega} W((\Rbb^h)^T ((\nabla \yb^h) \circ \Thetab^h\circ P^h)) dx \\
\nonumber & \geq & \liminf_{h \to 0} \Big[ \frac{1}{2} \int_\Omega Q_3(\chi_h \Gbb^h) dx + \frac{1}{h^8} \int_\Omega (1-\chi_h) W((\nabla \yb^h) \circ \Thetab^h\circ P^h)) dx \Big] \\
& \geq & \frac{1}{2} \int_{\Omega} \int_Y Q_3 (\Gbb(x,y))dy dx \geq \frac{1}{2} \int_{\Omega} \int_Y Q_2 (\Gbb''(x,y))dy dx.
\end{eqnarray}
Now by (\ref{zakor}) we have
\begin{eqnarray} \nonumber
\int_{-1/2}^{1/2} \int_Y Q_2 (\Gbb''(x_1,x_2,x_3,y))dy dx_3 &=& \int_YQ_2 (\Gbb_0(x_1,x_2,y))dy\\ & &\hspace{-5ex}+ \frac{1}{12} \int_Y Q_2 (\Gbb_1 (x_1,x_2,y))dy.
\end{eqnarray}
This implies the claim of the corollary.
\end{prooof}

Let us by $Q_2^H: \ZR^{2 \times 2}_{\sym} \times \ZR^{2 \times
2}_{\sym} \to \ZR_0^+$ denote the functional
\begin{equation} \label{defQ2h}
Q_2^H (\Gbb,\Fbb)=\min_{\small{\begin{array}{l}\ub_1 \in \dot{H}_{\#}^1(Y),\\ v \in \dot{H}_{\#}^2(Y) \end{array}}} I^H_{\Gbb,\Fbb},
\end{equation}
where we have by $I^H_{\Gbb,\Fbb}: (\dot{H}_{\#}^1(Y))^2 \times
\dot{H}_{\#}^2(Y) \to \ZR^+_0$ denoted the functional
\begin{eqnarray}\nonumber
I^H_{\Gbb,\Fbb}(\ub_1,v_1)&=&\int_Y \Big(Q_2(\Gbb+\Fbb \theta_0+\sym\nabla \ub_1(y)-\nabla^2
v_1(y)\theta_0(y))\\ & & \label{defIHGF}\hspace{5ex}+\frac{1}{12}Q_2(\nabla^2 v_1(y)) \Big)dy
\end{eqnarray}
By $I^L_0:W^{1,2}(\omega;\ZR^2) \times W^{2,2} (\omega) \to
\ZR^+_0$ we denote the functional
\begin{equation} \label{defil0}
I^L_0(\ub,v)=\frac{1}{2}\int_{\omega} Q_2^H(\sym \nabla \ub+\frac{1}{2} \nabla v \otimes \nabla v, -\nabla^2 v) dx+\frac{1}{24}\int_{\omega} Q_2(\nabla^2 v)dx.
\end{equation}
In the sequel we shall analyze the property of $Q_2^H$.
\begin{lemma} \label{ominimumu}
For every $\Gbb,\Fbb \in \ZR^{2 \times 2}$ there exists a
unique $\ub_1(\Gbb,\Fbb) \in (\dot{H}_{\#}^1(Y))^2$,
$v_1(\Gbb,\Fbb) \in \dot{H}_{\#}^2(Y)$ which minimizes the
functional $I^H_{\Gbb,\Fbb}$. This minimizer satisfies the
estimate
\begin{equation} \label{nnnnnnnnnn}
\|\ub_1(\Gbb,\Fbb)\|^2_{\dot{H}^1_{\#}(Y)}, \|v_1(\Gbb,\Fbb)\|^2_{\dot{H}^2_{\#}(Y)} \leq C(\|\Gbb\|^2+\|\Fbb\|^2),
\end{equation}
for some $C>0$. The functional $Q_2^H$ is a nonnegative
quadratic form. Let us by $V$ denote the subspace of
$\ZR^{2\times 2}_{\sym}$
\begin{eqnarray}
\nonumber V=\Big\{ \Abb= \left(\begin{array}{cc} a_{11} & a_{12} \\ a_{12} & a_{22} \end{array} \right); \forall y \in Y,  a_{11}\partial_{22} \theta(y)+a_{22} \partial_{11} \theta(y)-2a_{12} \partial_{12} \theta(y)=0 \Big\},
\end{eqnarray}
and by $V^{\perp}$ its orthogonal complement in $\ZR^{2 \times
2}_{\sym}$. Then there exists
 $C>0$ such that
\begin{equation}
Q_2^H(\Gbb,\Fbb)> C(\|\Gbb\|^2+\|\Fbb^{\perp}\|^2), \forall \Gbb, \Fbb \in \ZR^{2 \times 2}_{\sym},
\end{equation}
where $\Fbb^{\perp} \in V^{\perp}$ is the unique matrix such
that $\Fbb-\Fbb^{\perp} \in V$.
\end{lemma}
\begin{prooof}
The existence of the minimum of the functional
$I^H_{\Gbb,\Fbb}$ in the space $(\dot{H}_{\#}^1(Y))^2\times
\dot{H}_{\#}^2(Y)$ is guaranteed by the following facts:
\begin{enumerate}
\item for each $\Gbb,\Fbb$ the functional $I^H_{\Gbb,\Fbb}$
    is sequentially weakly lower semi continuous by the
    convexity of the form $Q_2$.
\item for each $\Gbb,\Fbb$ the functional $I^H_{\Gbb,\Fbb}$
    is coercive in the sense
    $$ \| \ub_1^n\|_{(\dot{H}_{\#}^1(Y))^2} \to +\infty \textrm{
    or } \|v_1^n\|_{\dot{H}_{\#}^2(Y)} \to +\infty \Rightarrow
I^H_{\Gbb,\Fbb}(\ub_1^n,v_1^n) \to +\infty. $$ To prove
 this let us note that the boundedness of
 $I^H_{\Gbb,\Fbb}(\ub_1^n,v_1^n)$ directly implies the
 boundedness of $\|v_1^n\|_{\dot{H}_{\#}^2(Y)}$. This, on
 the other hand, implies the boundedness of $\|\sym \nabla
 \ub_1^n\|_{L^2(Y)}$. From Korn's inequality we have the
 claim.
\end{enumerate}
The uniqueness of the minimum is the consequence of the strict
convexity of the form $Q_2$ on symmetric matrices.

 Let us by $\mathcal{A}:\ZR^{2\times2}_{\sym} \to
\ZR^{2\times2}_{\sym}$ denote the positive definite linear
operator which realizes the quadratic functional $Q_2$ i.e. we
have
\begin{equation}
Q_2(\Gbb)= (\mathcal{A} \Gbb,\Gbb), \forall \Gbb \in \ZR^{2 \times 2}_{\sym},
\end{equation}
where we have by $(\cdot,\cdot)$ denoted the standard scalar
multiplication on $\ZR^{2 \times 2}$. The minimization
formulation (\ref{defQ2h}) implies
\begin{eqnarray}\nonumber
& &\int_Y (\mathcal{A} (\sym \nabla \ub_1(y) -\nabla^2
v_1(y)\theta_0(y)), \sym \nabla \vphib(y))dy \\ \nonumber & &\hspace{10ex}+ \int_Y (\mathcal{A} (\sym \nabla \ub_1(y) -\nabla^2
v_1(y)\theta_0(y)), -\nabla^2 \upsilon(y)\theta_0(y))dy\\ \nonumber & & \hspace{20ex}+\frac{1}{6} \int_Y (\mathcal{A} \nabla^2 v_1(y), \nabla^2 \upsilon(y))dy\\ \nonumber
& &= -\int_Y (\mathcal{A}\Gbb,\nabla^2 \upsilon(y)\theta_0(y))dy+\int_Y (\mathcal{A} \Fbb \theta_0(y), \sym\nabla \vphib(y)-\nabla^2 \upsilon(y) \theta_0(y))dy \\ & & \label{nemaga} \hspace{40ex} \forall \vphib \in (\dot{H}_{\#}^1(Y))^2, \upsilon \in \dot{H}_{\#}^2(Y)
\end{eqnarray}
From (\ref{nemaga}) it can be easily seen that
$\ub_1(\Gbb,\Fbb)$, $v_1(\Gbb,\Fbb)$ depends linearly on
$(\Gbb,\Fbb)$ and thus $Q_2^H$ is a nonnegative quadratic form.
To see the estimate (\ref{nnnnnnnnnn}) we just take
$\vphib=\ub_1$, $\upsilon=v_1$. We obtain
\begin{equation} \label{nnnnn}
\| \ub_1-\nabla^2v_1 \theta_0 \|_{L^2(Y)}+\|\nabla^2 v_1\|_{L^2(Y)} \leq C(\|\Gbb\|+\|\Fbb\|).
\end{equation}
From (\ref{nnnnn}) we have (\ref{nnnnnnnnnn}). To check the
positive definiteness of $Q_2^H$ we have to check
\begin{equation}
Q_2^H(\Gbb,\Fbb)=0, \Gbb,\Fbb \in \ZR^{2 \times 2}_{\sym} \Rightarrow \Gbb=0, \ \Fbb \in V.
\end{equation}
Let us suppose $Q_2^H(\Gbb,\Fbb)=0$. From (\ref{defIHGF}) we
conclude that
\begin{equation} \label{jednaddddd}
v_1(y)=0, \quad \Gbb+\Fbb\theta_0(y) +\sym\nabla \ub_1(y)=0, \ \forall y \in Y.
\end{equation}
Integrating (\ref{jednaddddd}) in $Y$ and using the fact
$\ubb_1 \in \dot{H}^1_{\#}(Y)$ we conclude $\Gbb=0$.

From the second equation we conclude that $\Fbb \theta_0$ is
symmetrized gradient and this implies
\begin{equation}
\Fbb_{11} \partial_{22} \theta+\Fbb_{22} \partial_{11}\theta-2 \Fbb_{12} \partial_{12}\theta=0.
\end{equation}
From this  we have that $\Fbb \in V$.
\end{prooof}

Now we shall say something about the regularity of
(\ref{nemaga}). First we will need one technical lemma which is
a variant of well known lemma (see e.g. \cite{evans}) for
torus.
\begin{lemma} \label{pomocnaperiodicna}
Let us take $u \in L^2(Y)$ and let us extend $u$ to $\ZR^2$ by
periodicity.  Define the $i$-th difference quotient of size $h$
by
\begin{equation}
D^h_i(u)(y)=\frac{u(y+h\eb_i)-u(y)}{h}, \ (i=1,2)
\end{equation}
for $y \in Y$ and $h \in \ZR$ and let us define
\begin{equation}
D^h(u):=(D_1^h u,  D_2^h u).
\end{equation}
We have:
\begin{enumerate}[i)]
\item Suppose $u \in H_{\#}^1(Y)$. Then
$$ \|D^h u\|_{L^2(Y)}  \leq C\|Du\|_{L^2(Y)}, $$
for some $C>0$ and all $0<|h|<h_0$.
\item Assume $u \in L^2(Y)$ and there exists a constant $C$
    such that
    $$ \|D^h u\|_{L^2(Y)} \leq C,$$
for all $0<|h|<h_0$. Then
$$ u \in H_{\#}^1(Y), \ \textrm{with } \|Du\|_{L^2(Y)}
\leq C. $$
\end{enumerate}
\end{lemma}
\begin{prooof}
Let us take $Y'=[-1,2]^2$.
\begin{enumerate}[i)]
\item From standard theorem on difference quotients we
    conclude that there exists $C_1>$ such that
    \begin{equation}
 \|D^h u\|_{L^2(Y)}  \leq C_1\|Du\|_{L^2(Y')},
    \end{equation}
for all $0<|h|<1$. Since, by periodicity,
$\|Du\|_{L^2(Y')}=9\|Du\|_{L^2(Y)}$ we have the claim.
\item Let us take $u \in L^2(Y)$ and extend it, by
    periodicity to $L^2(Y')$. To conclude that
    $\|Du\|_{L^2(Y)} \leq C$ is the direct consequence of
    the standard theorem on difference quotient. We have to
    prove $u \in H_{\#}^1(Y)$, i.e. it has periodic
    boundary conditions. But this can be concluded from the
    fact that $\|D^h u\|_{L^2(V)} \leq 9C$, for every $V$
    open $Y \subset V \subset Y'$. This implies
    $\|Du\|_{L^2(V)} \leq 9C$, by standard theorem on
    difference quotients, and thus $u \in H^{1}(V)$. This,
    on the other hand, implies that $u$ has periodic
    boundary conditions.
\end{enumerate}
\end{prooof}
\begin{lemma} \label{regggg}
Let $\theta \in C^k_{\#}(Y)$. The solution $\ub_1, v_1$ of
(\ref{nemaga}) is a linear function of $\Gbb,\Fbb$ and $\ub_1$
is in the space $\dot{H}^{k+1}_{\#}(Y)$ and $v_1$ is in the
space $\dot{H}^{k+2}_{\#}(Y)$.
\end{lemma}
\begin{prooof}
Since we are on torus,  proving  regularity is easier since we
do not have boundary. Let us prove the claim for $k=1$. Let us
take $i \in \{1,2\}$ and test functions
$$ \vphib_{\alpha}=-D_{i}^{-h}(D_i^h \ub_\alpha), \quad
\alpha=1,2; \ \upsilon=-D_i^{-h}(D_i^h v_1)$$ in (\ref{nemaga}). From standard properties on difference quotient we conclude
\begin{eqnarray}\nonumber
& &\int_Y (\mathcal{A} (\sym \nabla D_i^h\ub_1(y) -\nabla^2
D_i^h v_1(y)\theta_0(y)), \sym \nabla D_i^h \ub_1)dy \\ \nonumber & &\hspace{10ex}+ \int_Y (\mathcal{A} (\sym \nabla D_i^h\ub_1(y) -\nabla^2
D_i^h v_1(y)\theta_0(y)), -\nabla^2 D_i^h v_1(y)\theta_0(y))dy\\ \nonumber & & \hspace{20ex}+\frac{1}{6} \int_Y (\mathcal{A} \nabla^2 D_i^h v_1(y), \nabla^2 D_i^h v_1(y))dy\\ \nonumber
& &= -\int_Y (\mathcal{A}\Gbb,\nabla^2 D_i^h v_1(y)\theta_0(y))dy \\ \nonumber & &\hspace{5ex}+\int_Y (\mathcal{A} \Fbb \theta_0(y), \sym \nabla D_i^h \ub_1(y)-\nabla^2 D_i^h v_1(y) \theta_0(y))dy \\
\nonumber& & \hspace{5ex} -\int_Y (\mathcal{A} (\sym \nabla \ub_1(y) -\nabla^2 v_1(y)D_i^h\theta_0(y)), \sym \nabla \ub_1(y))dy \\
\nonumber & & \hspace{5ex} -\int_Y (\mathcal{A} (\sym \nabla \ub_1(y) -\nabla^2
 v_1(y)\theta_0(y)), -\nabla^2  v_1(y)D_i^h\theta_0(y))dy\\
\label{zzzadnje} & & \hspace{5ex} -\int_Y (\mathcal{A} (\sym \nabla \ub_1(y) -\nabla^2
 v_1(y)D_i^h\theta_0(y)), -\nabla^2  v_1(y)\theta_0(y))dy
\end{eqnarray}
From (\ref{zzzadnje}), using the positive definiteness of
$\mathcal{A}$ and transposing the different quotient in the
first two terms on the right hand side, we conclude that there
exists $C>0$, dependent only on $\mathcal{A}$ such that
\begin{eqnarray} \nonumber
& &\|\sym \nabla D_i^h\ub_1 -\nabla^2
D_i^h v_1\theta_0) \|^2_{L^2(Y)}+ \|\nabla^2 D_i^h v_1\|^2_{L^2(Y)} \\ & & \leq C(\|\Gbb\|^2+\|\Fbb\|^2+\|\theta_0\|^2_{C^1(Y)})(\|\ub_1\|^2_{\dot{H}_{\#}^1(Y)}+\|v_1\|^2_{\dot{H}^2_{\#}(Y)}).
\end{eqnarray}
From this we conclude that $\|D_i^h \ub_1\|_{H^1(Y)}$, $\|D^h_i
v_1\|_{H^2(Y)}$ are bounded, independent of $h$, for $i \in
\{1,2\}$. Using lema \ref{pomocnaperiodicna} we have the claim
for $k=1$. For general $k$ we just differentiate (\ref{nemaga})
and repeat the arguments (see e.g. \cite{evans}). Thus we have
that the solution of (\ref{nemaga}) for $\theta \in
C^k_{\#}(Y)$ is given by
\begin{equation}
\ub_1^{\alpha}(y)=\Abb^\alpha_u(y) \cdot \Gbb+\Bbbb^\alpha_u(y) \cdot \Fbb, \ \alpha=1,2,\ v_1(y)=\Abb_v(y) \cdot \Gbb+\Bbbb_v(y) \cdot \Fbb,
\end{equation}
where $\Abb^{\alpha}_u,\Bbbb_u^{\alpha} \in
H^{k+1}_{\#}(Y;\ZR^{2 \times 2})$, $\Abb_v,\Bbbb_v \in
H^{k+2}_{\#}(Y;\ZR^{2 \times 2})$.
\end{prooof}

\subsection{Upper bound}
\begin{theorem}[optimality of lower bound] \label{upperbound}
Let  $v \in W^{2,2}(\omega)$,  $\ub \in W^{1,2} (\omega;
\ZR^2)$ and let $\theta \in C^2_{\#}(Y)$. Then for each
sequence $h \to 0$ there exists a subsequence, still denoted by
$h$ and appropriate $\yb^h \in W^{1,2}(\hat{\Omega}^h;\ZR^3)$
and $\Rbb^h \in W^{1,2}(\omega;\SO(3))$ such that
\begin{eqnarray}  \label{devide}
& &\| (\nabla \yb^h) \circ \Thetab^h\circ P^h -\Rbb^h \|_{L^2(\Omega)} \leq Ch^4, \\ \label{devide3}
& &\| \nabla \Rbb^h \|_{L^2(\omega)} \leq Ch^{2}, \frac{\Rbb^h-\Ibb}{h^2} \rightharpoonup \Abb \quad \textrm{in } W^{1,2}(\omega;\ZR^{3 \times 3}),
\end{eqnarray}
where $\Abb$ is defined by (\ref{defabb})
 and for $\ub^h$, $v^h$ defined by (\ref{definicijauv}) (where $\widetilde{\yb}$ should be replaced by $\yb$)  convergence (\ref{konvergencija})-(\ref{konvergencija2}) are valid
and
\begin{equation}
\lim_{h \to 0}  I^h (\yb^h)= I^L_{0}(\ub,v).
\end{equation}
\end{theorem}
\begin{prooof}
Let us assume that $\ub \in C^{\infty}(\omega;\ZR^2),v\in
C^{\infty}(\omega)$. Let us take an arbitrary $\ub_1 \in
\mathcal{D}(\omega;C^{\infty}_{\#}(Y))^2$ , $v_1
\in\mathcal{D}(\omega;C^{\infty}_{\#}(Y))$. Then we define
\begin{eqnarray} \nonumber
\yb^h (\Thetab^h(x_1,x_2,x_3^h)) &=& 
\Thetab^h (x_1,x_2,x_3)\\ \nonumber& &+\left(\begin{array}{c} h^4 \ub (x_1,x_2)+h^5 \ub_1(x_1,x_2,\per,x_3) \\ h^2 v(x_1,x_2)+h^4 v_1(x_1,x_2,\per,x_3) \end{array} \right)\\
\nonumber& &-h^2 x_3^h \left(\begin{array}{c} \partial_1 v(x')+h\partial_{y_1} v_1(x_1,x_2,\per,x_3)\\
\partial_2 v(x')+h\partial_{y_2} v_1(x_1,x_2,\per,x_3)\\ 0 \end{array} \right)\\
& &\nonumber\hspace{-5ex}-h^4\Big(\partial_1 v(x_1,x_2)+h\partial_{y_1} v_1(x_1,x_2,\per,x_3)\Big) \theta_0(\per) \eb_1\\
& &\nonumber\hspace{-5ex}-h^4\Big(\partial_2 v(x_1,x_2)+h\partial_{y_2} v_1(x_1,x_2,\per,x_3)\Big) \theta_0(\per) \eb_2 \\
\nonumber & &\hspace{-5ex}-h^3 x_3^h\Big(\partial_1 v(x_1,x_2) \partial_1 \theta(\per)+\partial_2 v(x_1,x_2) \partial_2 \theta(\per)\Big)\eb_3\\
 & &\hspace{-5ex}+ h^4x_3^h \db^0(x_1,x_2,\per)+\frac{1}{2} h^2 (x_3^h)^2 \db^1 (x_1,x_2,\per),
 \label{defhaty}
\end{eqnarray}
where $\db_0,\db_1 \in
\mathcal{D}(\omega;C^{\infty}_{\#}(Y))^3$ are going to be
chosen later.
 We calculate
\begin{eqnarray} \nonumber
\nabla \yb^h \nabla \Thetab^h &=&\nabla \Thetab^h + h^2 \Abb'+h^4\Bbbb'-h^2x_3^h \Cbb'\\ \nonumber & &
-h^3 \left( \begin{array}{ccc} \partial_1 v \partial_1 \theta & \partial_1 v\partial_2 \theta & 0 \\ \partial_2 v \partial_1 \theta  & \partial_2 v \partial_2 \theta & 0\\ 0 & 0 & \partial_1 v \partial_1 \theta+\partial_2 v \partial_2 \theta
  \end{array} \right)
  \\
   \nonumber & &
+h^4 \left( \begin{array}{ccc}- \partial_{y_1} v \partial_1 \theta & -\partial_{y_1} v\partial_2 \theta & 0 \\ -\partial_{y_2} v \partial_1 \theta  & -\partial_{y_2} v \partial_2 \theta & 0\\ \partial_{x_1}v_1 & \partial_{x_2}v_1 & 0
  \end{array} \right) \\
  \nonumber  & &- h^2 x_3^h \left( \begin{array}{ccc}  0 & 0 &0  \\
0 & 0 & 0 \\ \partial_1v\partial_{11} \theta+\partial_2v \partial_{12}\theta& \partial_1v\partial_{12} \theta+\partial_2v \partial_{22}\theta  & 0
  \end{array} \right)\\ \label{nablayt} & &
  + h^4 \db_0 \otimes \eb_3 +  h^2 x_3 ^h \db_1 \otimes \eb_3   +O(h^5),
\end{eqnarray}
where $\|O(h^5)\|_{L^\infty(\Omega)} \leq Ch^5$ and
\begin{eqnarray}
\Abb'&=&\left( \begin{array}{ccc} 0 & 0&-\partial_1 v-h\partial_{y_1}v_1 \\ 0& 0& -\partial_2 v-h\partial_{y_2}v_1 \\
   \partial_1 v+h\partial_{y_1}v_1  & \partial_2 v+h\partial_{y_2}v_1   & 0
  \end{array} \right), \\
\Bbbb' &=& \nonumber \left( \begin{array}{ccc} \partial_1  \ub^1+\partial_{y_1} \ub_1^1 &  \partial_2  \ub^1+\partial_{y_2} \ub_1^1   & 0 \\
   \partial_1  \ub^2+\partial_{y_1} \ub_1^2 & \partial_2  \ub^2+\partial_{y_2} \ub_1^2 & 0\\0 & 0 &0
  \end{array} \right) \\ & &-
  \left( \begin{array}{ccc} (\partial_{11} v+\partial_{y_1 y_1} v_1)\theta_0 & (\partial_{12} v+\partial_{y_1 y_2} v_1)\theta_0  & 0 \\
  (\partial_{12} v+\partial_{y_1 y_2} v_1)\theta_0 &   (\partial_{22} v+\partial_{y_2 y_2} v_1)\theta_0 & 0\\0 & 0 &0
  \end{array} \right), \\
\Cbb' &=& \left( \begin{array}{ccc}  \partial_{11}v+\partial_{y_1 y_1} v_1 & \partial_{12}v+\partial_{y_1 y_2} v_1 &0  \\
\partial_{12}v+\partial_{y_1 y_2} v_1 & \partial_{22}v+\partial_{y_2 y_2} v_1 & 0 \\ 0 & 0 & 0
  \end{array} \right).
\end{eqnarray}
 From (\ref{nablayt}), by using (\ref{inverz}), we
conclude
\begin{eqnarray} \nonumber
\nabla \yb^h &=& \Ibb+h^2\Abb'+ h^4 \Bbbb'-h^2 x_3^h\Cbb'\\ & &
 +h^4 \Obb
+h^2 x_3^h \Pbb
  + h^4 \db_0 \otimes \eb_3 + h^2 x_3 ^h \db_1 \otimes \eb_3  \label{nablayyt} +O(h^5),
\end{eqnarray}
where
\begin{equation}
\Obb = \left( \begin{array}{ccc}   0 & 0 & o^h_{13}  \\0 & 0 & o^h_{23} \\ o^h_{31} & o^h_{32}  & o_{33}^h
  \end{array} \right), \quad \Pbb = \left( \begin{array}{ccc} 0 &0 &0 \\ 0 &0 &0\\ p^h_{31} & p^h_{32} &0 \end{array} \right) \end{equation}
$o^h_{i3}, o^h_{3i}, p^h_{3 \alpha} \in
\mathcal{D}(\omega;C^{\infty}_{\#}(Y))$,
$\|o^h_{i3}\|_{L^\infty(\omega)},
\|o^h_{3i}\|_{L^\infty(\omega)}, \|p^h_{3
\alpha}\|_{L^\infty(\omega)} \leq C$, $C$ independent of $h$.

Let us define $$\Rbb^h:=e^{h^2 \Abb}.$$ The claims
(\ref{devide})-(\ref{devide3}) are easily checked to be valid
as well as the convergence (\ref{konvergencija}),
(\ref{konvergencija2}).

Using the identities $(\Ibb+\Abb)^T (\Ibb+\Abb)=\Ibb+2 \sym
\Abb +\Abb^T \Abb$ and $(\eb_3 \otimes \ab'-\ab' \otimes
\eb_3)^T (\eb_3 \otimes \ab'-\ab' \otimes \eb_3)=\ab' \otimes
\ab' +|\ab'|^2 \eb_3 \otimes \eb_3$ for $\ab' \in \ZR^2$ we
obtain
\begin{eqnarray} \nonumber
(\nabla \yb^h)^T (\nabla \yb^h) &=& \Ibb+ h^4[ 2\sym \Bbbb+ \nabla v \otimes \nabla v +|\nabla  v|^2 \eb_3 \otimes \eb_3]\\ \nonumber& & -2h^2 x_3^h (\nabla^2 v+\nabla^2_y v_1)+
2 h^4 \sym (\db_0 \otimes \eb_3)+\\ \nonumber & &+2h^2 x_3 ^h \sym(\db_1 \otimes \eb_3)
+2h^4\sym \Obb
+2h^2 x_3^h \sym \Pbb  +O(h^5).\\  \label{matereti} & &
\end{eqnarray}
For a symmetric $2 \times 2 $ matrix $\Abb''$ let
$c=\mathcal{L} \Abb'' \in \ZR^3$ denote the (unique) vector
which realizes the minimum in the definition of  $Q_2$. i.e.
$$ Q_2 (\Abb'')=Q_3 (\Abb''+\cb \otimes \eb_3+\eb_3 \otimes \cb). $$
Since $Q_3$ is positive definite on symmetric matrices, $\cb$
is uniquely determined and the map $\mathcal{L}$ is linear. We
choose now
\begin{eqnarray*}
\db_0^i &=& -\frac{1}{2} | \nabla v|^2 \delta_{i3}-\frac{1}{2} (\Obb_{i3}+\Obb_{3i}) + \mathcal{L} (\sym \Bbbb+ \frac{1}{2}\nabla  v \otimes \nabla  v), \\
\db_1^i &=& -\frac{1}{2} (\Pbb_{i3}+\Pbb_{3i}) -\mathcal{L} (\nabla^2 v+ \nabla^2_y v_1),
\end{eqnarray*}
where $\delta$ is Kronecker symbol, and we can see that all the
calculations are still valid. Taking the square root in
(\ref{matereti}) and using the frame indifference of $W$ and
the Taylor expansion we get
\begin{eqnarray}
\frac{1}{h^8}\int_{\hat{\Omega}^h} W(\nabla \yb^h)=\frac{1}{h^8}\int_{\hat{\Omega}^h} W([(\nabla \yb^h)^T \nabla \yb^h ]^{1/2}) \to I^L(\ub,v,\ub_1,v_1),
\end{eqnarray}

In the case $\theta \in C^{\infty}_{\#}(Y)$ we could use Lemma
(\ref{regggg}). Since we supposed only $\theta \in C^2_{\#}(Y)$
we continue as follows. Let us now take an arbitrary $\ub \in
W^{1,2}(\omega)$, $v \in W^{2,2}(\omega)$. Let $\ub_1(\nabla
\ub,-\nabla^2 v) \in L^2(\omega;\dot{H}^1_{\#}(Y))^2$, $v_1
(\nabla \ub,-\nabla^2 v) \in L^2(\omega;\dot{H}^2_{\#}(Y))$ be
from Lemma \ref{ominimumu}. We know take $\ub^n,v^n
,\ub_1^n,v_1^n$ smooth such that
\begin{eqnarray} \label{konuuuuu}
\|\ub^n-\ub\|_{W^{1,2}(\omega;\ZR^2)} &\leq& \frac{1}{n} \\ \label{konvvvvv}
\|v^n-v\|_{W^{2,2}(\omega)} &\leq& \frac{1}{n} \\
\| \ub_1^n-\ub_1 \|_{(L^2(\omega; \dot{H}_{\#}^1(Y))^2} &\leq& \frac{1}{n} \\
\| v_1^n-v_1 \|_{L^2(\omega; \dot{H}_{\#}^2(Y))} &\leq& \frac{1}{n} \\ \label{ooolll}
|I^L(\ub^n,v^n,\ub_1^n,v_1^n)-I^L_0(\ub,v)| &\leq& \frac{1}{n}
\end{eqnarray}
Now, we choose $\yb^{h_n} \in
W^{1,2}(\hat{\Omega}^{h_n};\ZR^3)$, from the proof, such that
\begin{enumerate}[i)]
\item $$ \|(\nabla\yb^{h_n}) \circ \Thetab^{h_n} \circ
    P^{h_n} -\Ibb\|_{L^2(\Omega)} \leq \frac{1}{n} $$
\item  for $v^{h_n}$ defined by (\ref{definicijauv}) we
    have $\|v^{h_n}-v^n\|_{W^{1,2}(\omega)} \leq
    \frac{1}{n}$.
\item  for $\tilde{\ub}^{h_n}$ defined by
\begin{eqnarray} & &\nonumber \tilde{\ub}^{h_n} (x_1,x_2):=\\ \nonumber &
&
\frac{1}{h_n^4}\Bigg(\int_{-1/2}^{1/2} \left( \begin{array}{c}{\yb}^{h_n}_1 \circ \Thetab^{h_n} \circ P^{h_n} \\ \yb^{h_n}_2 \circ \Thetab^{h_n} \circ P^{h_n} \end{array} \right)(x_1,x_2,x_3)dx_3-\left(\begin{array}{c} x_1 \\x_2\end{array}\right)\Bigg)\\
\nonumber& &+\left(\begin{array}{c}\partial_1 v^n(x_1,x_2) \theta_0(\per) \\ \partial_2 v^n(x_1,x_2)\theta_0(\per)\end{array}\right),
\end{eqnarray}
we have $\| \tilde{\ub}^{h_n}-\ub^n\|_{L^2(\omega;\ZR^2)}
\leq \frac{1}{n}$.
\item $$\|I^{h_n}(\yb^{h_n})-I^L(\ub^n,v^n,\ub_1^n,v_1^n)\|
    \leq \frac{1}{n}.$$
\end{enumerate}
By using ii) and iii) and compactness Lemma \ref{osnovna} and
Lemma \ref{ovavecer} we conclude that there exists $\Rbb^{h_n}$
such that (\ref{sranje1}), (\ref{sranje2}) is valid and for
$v^{h_n}$ and $\ub^{h_n}$ defined by (\ref{definicijauv}) we
have that (\ref{konvergencija})-(\ref{simetricno}) is valid (it
is easily seen from the properties ii) iii) and the assumptions
(\ref{konuuuuu}) and (\ref{konvvvvv}) that the limits of
$\ub^{h_n}$ and $v^{h_n}$ are $\ub$ and $v$). By using
(\ref{ooolll}) and iv) we conclude that $$ \lim_{n \to \infty}
I^{h_n}(\yb^{h_n})=I^L_0(\ub,v). $$ This finishes the proof of
theorem.
\end{prooof}

Lemma \ref{osnovna}, Lemma \ref{identifikacija} and Theorem
\ref{upperbound} enable us to standard theorem on convergence
of minimizers. We shall state it without proof (since it is
standard), assuming the external loads in $\eb_3$ direction.
For a more detailed discussion on external loads in the
standard
 F\"oppl-von K\'arm\'an case see \cite{Muller5}.

Let $\fb_3^h \in L^2(\hat{\Omega}^h;\ZR)$ be given with the property
\begin{eqnarray} \label{pretpostavkesile1}
\forall h, \frac{1}{h^6}\fb_3^h \circ \Thetab^h \circ P^h &=& \fb_3 \in L^2(\omega;\ZR),\\ \label{pretpostavkasile2} \int_{\hat{\Omega}^h} \fb_ 3^hdx&=&0.
\end{eqnarray}
It is not necessary to prove the equality in
(\ref{pretpostavkesile1}) neither that $\fb_3$ depends only on
$x_1,x_2$. For a more detailed discussion see the proof of
Theorem 2.5 in \cite{Lewicka1} (see also the proof of Theorem 6
in \cite{Velcic1}). Let us additionally assume
\begin{equation} \label{uvjetsile}
\int_{\omega} x_1 \fb_3 dx=0, \quad \int_{\omega} x_2 \fb_3 dx=0.
\end{equation}
 This can be
assumed by rotating the coordinate axis (the energy density $W$
in rotational invariant).
From (\ref{determinanta}) it can be concluded that
\begin{equation} \label{uvjetsile1}
|\int_{\hat{\Omega}^h} x_1 \fb_3^h dx| \leq Ch^4, \quad |\int_{\hat{\Omega}^h} x_2 \fb_3^h dx| \leq Ch^4.
\end{equation}
 The total energy functional $J^h$
(divided by $h^2$),  defined on the space
$W^{1,2}(\hat{\Omega}^h;\ZR^3)$, is given by
\begin{equation} \label{defjh} J^h (\yb^h)= h^8 I^h(\yb^h)-
\frac{1}{h^2} \int_{\hat{\Omega}^h} \fb_3^h \yb_3^h.
\end{equation}
The following theorem is the main result and its proof follows
the proof of Theorem 2 in \cite{Muller3}.
\begin{theorem}[$\Gamma$-convergence] \label{najglavnijiteorem}
 Let us suppose that $\theta \in C^2_{\#}(Y)$ and
  $\fb_3^h \in L^2(\hat{\Omega}^h;\ZR)$ is given and satisfies (\ref{pretpostavkesile1}), (\ref{pretpostavkasile2}) and (\ref{uvjetsile}). Then:
\begin{enumerate}[1.]
\item  There exists $C_1,C_2>0$ such that for every $h>0$
    we have
\begin{equation} \label{nejinf} C_1 h^2 \geq \inf\left\{ \frac{1}{h^8} J^h (\yb^h) ; \yb^h \in W^{1,2} ( \hat{\Omega}^h;\ZR^3) \right\} \geq -C_2.   \end{equation}
In the case $\langle \theta \rangle \neq 0$ we can take
$C_1=0$.
\item If $\yb^h \in W^{1,2}(\hat{\Omega}^h;\ZR^3)$ is a
    minimizing sequence of $\frac{1}{h^8} J^h$, that is
\begin{equation} \label{svojmin} \lim_{h \to 0} \Big( \frac{1}{h^8} J^h (\yb^h)- \inf \frac{1}{h^8} J^h \Big)=0, \end{equation}
then we have that there exists $\bar{\Rbb}^h \in \SO(3), \ \cb^h \in \ZR$ such that  the sequence $(\bar{\Rbb}^h, \yb^h)$ has its subsequence (also not relabeled) with the following property:
\begin{enumerate}[i)]
\item For the sequence $
    \widetilde{\yb}^h:=(\bar{\Rbb}^h)^T \yb^h- \cb^h$
    and $\ub^h, v^h$ defined by (\ref{definicijauv})
    the following is valid
    \begin{eqnarray*}
\ub^h &\rightharpoonup& \ub \quad \textrm{weakly in } W^{1,2}(\omega;\ZR^2), \\
v^h &\to& v  \quad \textrm{in } W^{1,2}(\omega),\  v \in W^{2,2}(\omega).
\end{eqnarray*}
\end{enumerate}
 Any accumulation point $(\ub,v,\bar{\Rbb})$ of the
 sequence $(\ub^h,v^h,\bar{\Rbb}^h)$ minimizes the
 functional
\begin{equation}
J^L_0(\ub,v,\bar{\Rbb})=I^L_0(\ub,v)- \bar{\Rbb}_{33} \int_{\omega}  \fb_3(x_1,x_2) \big(v(x_1,x_2)+\theta(\per)\big)dx_1 dx_2,
\end{equation}
where $I^L_0$ is defined in (\ref{defil0}). Moreover, for
$\fb_3 \neq 0$ we have $\bar{\Rbb}_{33}=1$ or
$\bar{\Rbb}_{33}=-1$.
\item The minimum of the functional  $J^L_0$ exists in the
    space $W^{1,2}(\omega;\ZR^2) \times W^{2,2}(\omega)
    \times \SO(3)$. If $\yb^h \in
    W^{1,2}(\hat{\Omega}^h;\ZR^3)$ is a minimizing sequence
    (not relabeled) of $\frac{1}{h^8} J^h$ then we have
    that
\begin{equation}
    \lim_{h \to 0} \frac{1}{h^8} J^h (\yb^h)= \min_{\ub \in W^{1,2}(\omega;\ZR^2), \ v \in W^{2,2}(\omega),\ \bar{\Rbb} \in \SO(3)} J^L_0 (\ub,v, \bar{\Rbb}).
\end{equation}
\end{enumerate}
\end{theorem}
\begin{remark} \label{jjjakovazno}
When we compare this model with  the ordinary plate model of
the F\"oppl-von K\'arm\'an type we see that in the energy
expression we mix the term $\nabla^2 v$ which measures the
bending of the plate with the term $\sym \nabla \ub+\frac{1}{2}
\nabla v \otimes \nabla v$ which measures the stretching of the
plate. Thus the imperfect plate (periodically wrinkled) can
cause this kind of behavior.
\end{remark}

\begin{remark}
In the same way we could derive the linear model of
periodically wrinkled plate, starting from $3D$ nonlinear
theory (see \cite{Muller3}; see also \cite{Velcic1}). The
linear model would be like this model without the term $\nabla
v \otimes \nabla v$ in the expression (\ref{defil0}), where we
defined $I^L_0$.
\end{remark}

\end{document}